\documentstyle[amssymb,amsfonts,12pt]{amsart}

\newtheorem{theorem}{Theorem}[section]
\newtheorem{proposition}[theorem]{Proposition}
\newtheorem{lemma}[theorem]{Lemma}
\newtheorem{observation}[theorem]{Observation}
\newtheorem{corollary}[theorem]{Corollary}
\newtheorem{definition}[theorem]{Definition}

\def\C{{\mbox{\rm\kern.24em
\vrule width.03em height1.43ex depth-.052ex \kern-.26em C}}}
\def\Z{{\bf Z}}
\def\R{{\mbox{\rm I\kern-.22em R}}}

\def\Q{{\bf Q}}
\def\T{{\bf T}}
\def\I{{\bf I}}
\def\J{{\bf J}}

\def\size{{\rm size}}
\def\diam{{\rm diam}}
\def\sgn{{\rm sgn}}
\def\osc{{\rm osc}}
\def\BMO{{\rm BMO}}

\def\m{{\bf m}}
\def\K{J}  

\def\pv{{\vec P}}

\def\Pv{{\vec{\bf P}}}
\def\dist{{\rm dist}}

\def\111{\gamma}
\def\hyperplane{\Gamma}

\def\subspace{{\Gamma'}}

\def\be#1{\begin{equation}\label{#1}}
\def\bas{\begin{align*}}
\def\eas{\end{align*}}
\def\bi{\begin{itemize}}
\def\ei{\end{itemize}}
\newenvironment{proof}{\noindent {\bf Proof} }{\endprf\par}
\def \endprf{\hfill  {\vrule height6pt width6pt depth0pt}\medskip}
\def\emph#1{{\it #1}}

\title{Uniform estimates on multi-linear operators with
modulation symmetry}

\author{Camil Muscalu}
\address{Department of Mathematics, UCLA, Los Angeles CA 90095-1555}
\email{camil@@math.ucla.edu}

\author{Terence Tao}
\address{Department of Mathematics, UCLA, Los Angeles CA 90095-1555}
\email{tao@@math.ucla.edu}

\author{Christoph Thiele}
\address{Department of Mathematics, UCLA, Los Angeles CA 90095-1555}
\email{thiele@@math.ucla.edu}

\begin{document}

\begin{abstract}  In a previous paper \cite{mtt-1} in this series, we gave $L^p$ estimates for multi-linear operators given by multipliers which are singular on a non-degenerate subspace of some dimension $k$.  In this paper we give uniform estimates when the subspace approaches a degenerate region in the case $k=1$, and when all the exponents $p$ are between $2$ and $\infty$.  In particular we recover the non-endpoint uniform estimates for the Bilinear Hilbert transform in \cite{grafakos-li}.
\end{abstract}

\maketitle

\section{Introduction}\label{intro-sec}

We are concerned with $n$-linear forms mapping $n$ Schwartz functions
on the real line to a complex number. We shall assume these
forms are invariant under a simultaneous translation of all $n$
functions. Dually, these forms can be viewed as operators mapping
$n-1$ functions to a distribution. Our goal is to prove
$L^p$ regularity for these forms and operators.

By the Schwartz kernel theorem, we can identify such an $n$-linear
form with a distribution in $\R^n$. Translation invariance implies that
the Fourier transform of this distribution lives on the hyperplane
$$ \hyperplane := \{ \xi \in \R^n: \xi_1 + \ldots + \xi_n = 0 \}.$$
Since we are interested in $L^p$ regularity, we may restrict attention
to the case that this distribution is a function $m(\xi_1, \ldots, \xi_n)$ on this
hyperplane. 

Then the associated multi-linear form $\Lambda := \Lambda_m$ is given by
$$ \Lambda_m(f_1, \ldots, f_n) := \int \delta(\xi_1 + \ldots + \xi_n)
m(\xi) \hat f_1(\xi_1) \ldots \hat f_n(\xi_n)\ d\xi $$
where $\xi =: (\xi_1, \ldots, \xi_n)$ and $\delta$ denotes the Dirac delta.  
Dually, the
associated multi-linear operator $T := T_m$ is given by
\be{tm-def}
 T_m(f_1, \ldots, f_{n-1}\hat{)}(-\xi_n) := \int \delta(\xi_1 + \ldots + \xi_n) m(\xi) \hat f_1(\xi_1) \ldots \hat f_{n-1}(\xi_{n-1})\ d\xi_1 \ldots d\xi_{n-1},
\end{equation}
the relationship between $T$ and $\Lambda$ is given by
\be{lambda-t} \Lambda(f_1, \ldots, f_n) = \int T(f_1, \ldots, f_{n-1})(x) f_n(x)\ dx.
\end{equation}

Examples of such objects include the pointwise product operator 
$$ T_1(f_1, \ldots, f_{n-1}) := f_1 \ldots f_{n-1},$$
which occurs when the multiplier $m$ is identically 1.  
Operators from classical paraproduct theory studied in \cite{coifmanm1}-\cite{coifmanm6}, \cite{kenigs}, \cite{grafakost} 
also fall into this category and are given by
multipliers satisfying symbol estimates as in \eqref{symb-gamma-old}
below with $\Gamma'=\{0\}$. When $n=2$ these operators are just linear Fourier 
multipliers.  More recently, multipliers satisfying \eqref{symb-gamma-old}
with nontrivial subspaces $\Gamma'$ have been studied. 
Observe that \eqref{symb-gamma-old} is invariant under
translations in direction of $\Gamma'$, so the class of multipliers
satisfying this condition is translation invariant.
Taking the Fourier transform we obtain a class of forms and operators
which has modulation symmetry. In case the multiplier is invariant
under translations in direction of $\Gamma'$, the operator itself
has a modulation symmetry. Hence the title of this paper.
Such forms and operators have been discussed in 
\cite{mtt-1}, \cite{gilbert-nahmod-1} (also \cite{tao:xsb}, \cite{mtt-biest}, \cite{mtt-walsh-nest}, \cite{gilbert-nahmod-walsh}, \cite{gilbert-nahmod-2}).

Given a multiplier $m$, the main interest in the subject is to
obtain $L^p$ estimates for $\Lambda_m$ of the form
\be{lambda-bound}
| \Lambda_m(f_1, \ldots, f_n) | \leq C \prod_{i=1}^n \|f_i\|_{p_i}
\end{equation}
when $1 \leq p_i \leq \infty$.  
By duality this is equivalent to $T_m$ having the mapping property
$$ T_m: L^{p_1} \times \ldots \times L^{p_{n-1}} \to L^{p'_n}$$
on test functions, where $p'$ is given by $1/p + 1/p' := 1$.
When \eqref{lambda-bound} obtains, we say that $\Lambda_m$ is of \emph{strong type} $(1/p_1, \ldots, 1/p_n)$. In the dual formulation for $T_m$ one may
also consider the case ${p_n}'<1$, but we shall not do so in this paper.

In the pointwise product case $m \equiv 1$ one has strong type whenever $1 \leq p_i \leq \infty$ and one has the scaling condition
\be{scaling}
\sum_{i=1}^n \frac{1}{p_i} = 1.
\end{equation} 
These estimates (with the exception of some of the endpoint estimates) 
also generalize to paraproducts and the bilinear Hilbert transform.  We cite the following multiplier theorem, proven\footnote{Actually, a more general theorem was proven in \cite{mtt-1} in which some $1/p_i$ were allowed to be zero or negative, but we will not discuss this further here.  Also, the $n=3$ case of Theorem \ref{mtt-thm} was first proven in \cite{gilbert-nahmod-1}, \cite{gilbert-nahmod-2}.} in \cite{mtt-1}:

\begin{theorem}\label{mtt-thm}\cite{mtt-1}  Let $\subspace$ be a subspace of $\hyperplane$ of dimension $k$ where
\be{k-constraint}
0 \leq k < n/2. 
\end{equation}
Assume that $\subspace$ is non-degenerate in the sense that for every
$1 \leq i_1 < i_2 < \ldots i_k \leq n$, the space $\subspace$ is a graph
over the variables $\xi_{i_1}, \ldots, \xi_{i_k}$.
Suppose that $m$ satisfies the estimates
\be{symb-gamma-old}
|\partial_\xi^\alpha m(\xi)| \leq C \dist(\xi,\subspace)^{-|\alpha|}
\end{equation}
for all partial derivatives $\partial_\xi^\alpha$ on $hyperplane$ up to some sufficiently large
finite order.  Then \eqref{lambda-bound} holds whenever \eqref{scaling} holds and $1 < p_i < \infty$ for all $i = 1, \ldots, n$.
\end{theorem}

Following up on the work begun in \cite{thiele},
we consider in this paper the problem of obtaining uniform estimates when the subspace $\subspace$ becomes increasingly degenerate.  In order to do so one must modify the condition \eqref{symb-gamma-old}.  To illustrate this,
suppose that $n=3$, $k=1$, and $\subspace$ has degenerated completely to 
$$ \subspace = \{ (\xi_1,\xi_2,\xi_3) \in \hyperplane: \xi_1 = 0 \}.$$
(This space is a graph over the $\xi_2$ or $\xi_3$ variables, but not over $\xi_1$).  Let $m$ denote the multiplier
$$ m(\xi_1,\xi_2,\xi_3) := \psi(\xi_1) a(\xi_2)$$
where $\psi$ is a bump function on $[1,2]$ and $a$ is any smooth function satisfying the estimates $\| \partial_\xi^k a \|_\infty \leq C_k$ for all
integers $k$.  This multiplier $m$ satisfies $\eqref{symb-gamma-old}$
for the given degenerate subspace $\subspace$.
Then one has
$$ \Lambda_m(f_1, f_2, f_3) = \int T_\psi(f_1) T_a(f_2) f_3$$
where $T_m$ is the Fourier multiplier corresponding to $m$.  This operator can be quite badly behaved because $a$ need not be an $L^p$ multiplier for any $p \neq 2$.  Indeed, it is easy to construct examples which show unboundedness of this operator whenever $p_2 > 2 > p_3$ or $p_2 < 2 < p_3$.

The main result of this paper is the following $k=1$ result.

\begin{theorem}\label{main}  Let $n \geq 3$, and let $\subspace = span( v )$ be a one-dimensional subspace of $\hyperplane$, where $v = (v_1, \ldots, v_n)$ and $v_1, \ldots, v_n$ are non-zero real numbers which sum to zero.  Define the metric $d_v(x,y)$ on $\hyperplane$ by
$$ d_v(x,y) := \sup_{1 \leq i \leq n} \frac{|x_i - y_i|}{|v_i|}$$
and write
$$ d_v(x,\subspace) := \inf_{y \in \subspace} d_v(x,y).$$
Suppose that $m$ satisfies the estimates
\be{symb-gamma}
|\partial_\xi^\alpha m(\xi)| \leq C \prod_{i=1}^n 
(v_i d_v(\xi,\subspace))^{-\alpha_i}
\end{equation}
for all partial derivatives $\partial_\xi^\alpha$ on $\hyperplane$ up to some finite order.  Then \eqref{lambda-bound} holds whenever \eqref{scaling} holds and $2 < p_i < \infty$ for all $i = 1, \ldots, n$, with the bounds uniform in the choice
of $v_i$.
\end{theorem}

We do not know how to modify \eqref{symb-gamma} in the higher rank case $k > 1$, mainly because we have no natural analogue of the metric $d_v$.  In analogy with \cite{thiele}, \cite{grafakos-li-2} one expects that one should be able to go beyond the case $2 < p_i < \infty$, but we do not pursue these matters here. Restraining
ourselves to the case $2<p_i<\infty$ gives us some considerable technical
simplification.

In the non-degenerate case, when all the $v_i$ have comparable magnitudes, this theorem is a corollary of Theorem \ref{mtt-thm}.  A more careful examination of the proof of this Theorem in \cite{mtt-1} would reveal that the constant given for \eqref{lambda-bound} would grow polynomially in the ratio between the largest and smallest magnitudes of $v_i$.

A special case of this theorem occurs when $n=3$, $v = (\beta_2- \beta_3, \beta_3 - \beta_1, \beta_1 - \beta_2)$, for some distinct real numbers $\beta_1, \beta_2, \beta_3$,
and $m(\xi_1,\xi_2,\xi_3) = \sgn(\beta_1 \xi_1 + \beta_2 \xi_2 + \beta_3 \xi_3)$.  This corresponds to the (essentially one-parameter) family of bilinear Hilbert transforms
\be{bht}
\Lambda_m(f_1,f_2,f_3) = C \int \int \prod_{i=1}^3 f_i(x - \beta_i t) \frac{dt}{t} dx.
\end{equation}
In \cite{laceyt1}, \cite{laceyt2} the estimate \eqref{lambda-bound} was proven for all $1 < p_i < \infty$ obeying \eqref{scaling}, but with the bound growing polynomially when two of the $\beta_i$ approached each other.  In \cite{thiele} the weak-type estimate $T_m: L^2 \times L^2 \to L^{1,\infty}$ was shown with uniform control on the bounds when $\beta_1$ or $\beta_2$ approached $\beta_3$, but not when $\beta_1$ approached $\beta_2$.  (This weak-type bound is already sufficient to prove a version of Calder\'on's conjecture which is strong enough to recover the boundedness of the Cauchy integral on Lipschitz curves.  See \cite{thiele} for further discussion).  More recently, \eqref{lambda-bound} was established whenever $2 < p_i < \infty$ obeyed \eqref{scaling}, with bounds uniform over all values of $\beta$ in \cite{grafakos-li}.  Thus Theorem \ref{main} is already known in this case.   Also, estimates 
beyond the $2 \le p_i \le \infty$ case have been obtained 
in \cite{grafakos-li-2}.

Our argument follows the standard approach to the modulation
invariant setting of \cite{laceyt1}, \cite{laceyt2}, \cite{thiele}, \cite{mtt-1} (see also \cite{fefferman}, \cite{grafakos-li}, \cite{grafakos-li-2}).  We perform dyadic decompositions in both space and frequency, which has the effect of decomposing the phase plane into an overlapping set of tiles.  We then split $\Lambda_m$ into pieces associated to various $n$-tuples of tiles.  By using 
tree selection arguments and an orthogonality argument based on Bessel's inequality for tiles as in the above references, we can reduce matters to estimating the contribution of a single tree of tiles.  After eliminating the spatial cutoffs, the problem now becomes that of obtaining uniform estimates for paraproducts.  This estimate may be of independent interest, it has been shown in the
prequel \cite{mtt-2} of this paper.

We quote it in the form that is needed here:

\begin{proposition}\label{paraproduct-cor}  Let $n \geq 2$, and let $M_1, \ldots, M_n$ be real numbers.  For each $1 \leq i \leq n$ and $j \in \Z$, let $\pi_{j,i}$ be a Fourier multiplier whose symbol is a bump function adapted to $\{ \xi: |\xi| \le 2^{j + M_i} \}$.  Suppose that for each $k \in \Z$ there exists at least one $1 \leq i \leq n$ such that the symbol of $\pi_{j,i}$ vanishes at the origin. Then one has the estimate
\be{para-cor}
\sum_j |\int \prod_{i=1}^n \pi_{j,i} f_i|
\leq C_{(p_i),n} \prod_{i=1}^n \|f_i\|_{p_i}
\end{equation}
for all $1 < p_i < \infty$ obeying \eqref{scaling}, where the constant $C_{(p_i),n}$ 
depends on the $p_i$, and $n$ but is independent of the $M_i$ and $f_i$.
\end{proposition}

Here we say a bump function is adapted to an interval $I$ if it is
supported in this interval and its $k$-th derivative is bounded by
$|I|^{-k}$ for all $k$ up to some sufficiently large power, which
may depend on $n$ and $(p_i)$. 
Proposition \ref{paraproduct-cor} may be viewed as a lacunary 
version of Theorem \ref{main}, in the same way that paraproduct estimates are a lacunary version of the results in \cite{laceyt1}, \cite{laceyt2}, \cite{mtt-1} and the theory of maximal truncated Hilbert transforms are a lacunary version of the Carleson-Hunt theorem.  Indeed, it is possible to derive Proposition
\ref{paraproduct-cor} as a special case of Theorem \ref{main}, at least in the $2 < p_i < \infty$ case.  

One of the main innovations in this paper lie in using phase space projections to obtain the tree estimate from Proposition \ref{paraproduct-cor}. Such phase
plane projections were previously only known and utilized in the Walsh case
\cite{thiele2}. The fact that we are restricted to the $2 < p_i < \infty$ case allows for some simplifications in the argument; 
most notably the argument works with a uniform spatial localization
of all functions involved, unlike in \cite{thiele}. 
Morover, one can replace all the $f_i$ by characteristic functions.  
Also, since all functions are locally in $L^2$, one does not need to 
remove any exceptional set (which needs to be done for the $p_i < 2$ theory).

The first author was partially supported by a Sloan Dissertation Fellowship.
The second author is a Clay Prize Fellow and is supported by grants from
the Sloan and Packard Foundations. The third author was partially supported
by a Sloan Fellowship and by NSF grants DMS 9985572 and DMS 9970469.
The third author would like to thank the Mathematics Department
of Arizona State University for their hospitality during a visit
in which part of the work on this paper was done.

\section{Discretization in the frequency space}

The first step in the proof of Theorem \ref{main} is a Whitney
decomposition away from $\Gamma'$ of the multiplier $m$, which is 
defined as a function in frequency variables. The Whitney pieces
will be smooth functions and compactly supported in the frequency 
variables. While we would like to think of these pieces as parameterized
by dyadic intervals, it will not be possible to use the standard 
dyadic grid to parameterize them. As in \cite{laceyt1}, we will
spend some effort in constructing a grid structure akin to the 
dyadic structure. We thought it desirable to not having to do
this technical step, which involves a strong interaction of 
intervals in $i$- and $i'$- coordinates for $i\neq i'$, but 
we have not been able to proceed without it.

We fix $n$. All constants that follow (typically denoted by $C$)
may depend on $n$.
We also fix $p_1,\dots ,p_n$ as in Theorem \ref{main}, and
all constants are also allowed to depend on $p_1,\dots, p_n$.
In particular, we choose a large $N$ depending on these 
exponents, which will describe decay of functions in physical space.

We shall need the constants $C_0=2^{100n}$ and 
$\K=2^{C_0}$. These measure the fineness of our Whitney decomposition
and separation in scales respectively. All constants $C$ can be viewed
as dependent on $C_0$ and $\K$.

If $Q$ is a box in $\R^n$ with sides parallel to the axes, we use 
$Q_1, \ldots, Q_n$ to denote the intervals comprising $Q$, thus 
$Q = Q_1 \times \ldots \times Q_n$.  All our intervals and boxes 
will be closed, but when we say two intervals are disjoint we mean
disjoint up to possibly points at the boundary. 
For $A>0$, we denote by $AQ$ the box with the same center
as $Q$ and $A$ times the sidelength.

In the non-degenerate case treated in \cite{mtt-1}, we used the standard Euclidean Whitney decomposition of $\R^n\setminus \Gamma'$ 
into cubes.  This decomposition does not work well in the near-degenerate setting; following \cite{thiele} we shall instead decompose the multiplier into boxes which are adapted to the subspace $\Gamma'$.  To do this, it is convenient to perform a rescaling to convert the near-degenerate space $\subspace$ into a non-degenerate one, which we choose to be the diagonal
\[\tilde \Gamma' =\{(\xi,\ldots, \xi):\xi\in \R\}\subset \R^n\ \ \ .\]

We shall use tildes to denote quantities that are defined in the rescaled setting.

Fix $v_1, \ldots, v_n$.  By symmetry we may assume $|v_1|\ge \dots\ge |v_n|$.
By rescaling we may assume $|v_n|=1$. It is 
important that all our constants $C$ will be independent of the $v_i$; 
indeed, this is the whole point of this paper. We also define
$M_i$ and $\m_i$ by $2^{M_i}=|v_i|$ and $\K \m_i=M_i$.

Let $L: \R^n \to \R^n$ be the linear transformation
$$ L(x_1, \ldots, x_n) := (v_1 x_1, \ldots, v_n x_n).$$
The space 
$L^{-1}(\subspace)$ is thus the diagonal $\tilde{\Gamma}'$ in $\R^n$.

One can easily verify from \eqref{symb-gamma} that one can find a smooth 
function $\tilde m: \R^n \backslash \tilde \Gamma ' \to \R$ which agrees with $m \circ L$ on $L^{-1}(\hyperplane)$ and which satisfies the symbol estimates
$$ | \partial_\xi^\alpha \tilde m(\xi)| \leq C 
\dist(\xi,\tilde \Gamma ')^{-|\alpha|}$$
for all $\alpha$ up to some large finite order $N^5$, where $\dist$ is now Euclidean distance.  

Let $\tilde \Q$ denote the set of all cubes $\tilde Q$ which have side-length $2^{j}$ for some integer $j$, whose center lies in the lattice $2^{j-10} \Z^n$, and which obey the Whitney conditions
\begin{equation}\label{whitney1}
 C_0 \tilde Q\cap \tilde\Gamma '=\emptyset
\end{equation}
\begin{equation}\label{whitney2}
 4 C_0 \tilde Q\cap \tilde\Gamma '\neq \emptyset\ .
\end{equation}

The sets $\frac{1}{10}\tilde Q$ form a finitely overlapping cover of $\R^n \backslash \tilde \Gamma '$, and so one may decompose
$$ \tilde m = C \sum_{\tilde Q \in \tilde \Q} \tilde m_{\tilde Q}$$
where each $\tilde m_{\tilde Q}$ is a bump function adapted to $\frac{1}{2} \tilde Q$ with degree of regularity $N^4$.  
I.e., $\tilde m_{\tilde Q}$ is supported in $\frac 12 \tilde{Q}$
and satisfies
\[|\partial_\xi^\alpha \tilde m_{\tilde Q}(\xi)|\le \diam (\tilde Q)^{-|\alpha|}\]
for all $|\alpha|\le N^4$. It thus suffices to show that

$$
|\sum_{\tilde Q \in \tilde \Q} \int \delta(\xi_1 + \ldots + \xi_n)
\tilde m_{\tilde Q}(L^{-1}(\xi)) \hat f_1(\xi_1) \ldots \hat f_n(\xi_n)\ d\xi|
\leq C \prod_{i=1}^n \|f_i\|_{p_i}.$$

For each $\tilde Q \in \tilde \Q$, we may use a Fourier series and the smoothness of $\tilde m_{\tilde Q}$ to decompose $\tilde m_{\tilde Q}$ into tensor products
$$ \tilde m_{\tilde Q} = C 
\sum_{k \in \Z^n} (1+|k|)^{-10n} \prod_{i=1}^n \tilde m_{\tilde Q,i,k}(\xi_i)$$
where the $m_{\tilde Q,i,k}$ are bump functions adapted to the interval $\tilde Q_i$ with degree of regularity $N^3$ uniformly in $k$.  It thus suffices to show that
\begin{equation}\label{tilde-sum}
|\sum_{\tilde Q \in \tilde \Q} \int \delta(\xi_1 + \ldots + \xi_n)
\prod_{i=1}^n \tilde m_{\tilde Q,i,k}(v_i^{-1}\xi_i) \hat f_i(\xi_i)\ d\xi|
\leq C \prod_{i=1}^n \|f_i\|_{p_i}
\end{equation}
for each $k$. We fix $k$ once and for all and for notational convenience
drop the index $k$, i.e., write $\tilde m_{\tilde Q,i}$
instead of $\tilde m_{\tilde Q,i,k}$.

Next, we will make the collection of cubes sparser.
\begin{definition}
We call a collection $\tilde{\Q'}\subseteq \tilde{\Q}$ sparse, if 
we have for any $\tilde Q,\tilde Q'\in \tilde \Q'$ 
and any $1\le i\le n$ the following properties
\begin{equation}\label{sparse-scale}
|\tilde Q_i|< |\tilde Q_i'|\implies |\tilde Q_i|\le 2^\K 
|\tilde Q_i'|\ ,
\end{equation}
\begin{equation}\label{sparse-single-scale}
|\tilde Q_i|=|\tilde Q_i'|,\ 
\tilde Q_i\neq \tilde Q_i'
 \implies \dist(\tilde Q_i,\tilde Q_i')\ge 2^\K\ 
|\tilde{Q_i}|\ ,
\end{equation}
\begin{equation}\label{sparse-injectivity}
\tilde Q_i=\tilde Q_i' 
\implies \tilde Q=\tilde{Q}'\ .
\end{equation}
\end{definition}

Thanks to the fact that at given size the cubes in $\tilde{\Q}$
form essentially a one dimensional set in direction of the diagonal, 
we may decompose $\tilde \Q$ into a bounded number of sparse sets. Thus
it suffices to prove \eqref{tilde-sum} under the assumption that
$\tilde{\Q}$ is now a sparse subset of the original $\tilde{\Q}$.
For convenience of notation we shall call this sparse set again $\tilde \Q$.
By a standard limiting argument, we can also assume that $\tilde{\Q}$
is finite, as long as the estimates do not depend on the set $\tilde{\Q}$.

We now begin to introduce a structure akin to a dyadic grid.
We shall be interested in the enlarged cubes $1000\tilde Q$.
For each interval $1000 \tilde{Q}_i$ we shall
define a slightly (by at most one percent on either side)
increased  interval $\overline {Q}_i\supset 1000 \tilde{Q}_i$ 
so that these increased intervals have good dyadic properties
as formulated in Lemma \ref{dyadic-lemma}.

By \eqref{sparse-scale} the possible sizes of cubes in $\tilde\Q$
come in discrete quantities, each separated by a large factor. We 
shall refer to these different sizes as different scales.

We shall define the intervals $\overline{Q}_i$ 
successively beginning with the smallest scale in $\tilde \Q$.
For $\tilde Q\in \tilde \Q$ of the smallest scale we simply
set $\overline Q_i=1000\tilde Q_i$.

Observe that the distance between
$\tilde{Q}_i$ and $\tilde{Q}_j$ for $i\neq j$ is less than
$8C_0 |\tilde{Q}_i|$, because the cartesian product of
$4C_0 \tilde{Q_i}$ and $4C_0 \tilde{Q_j}$ contains a point
on the diagonal by \eqref{whitney2}.
Thus the diameter of the
convex hull $h(\overline{Q})$ of all $10\overline Q_i$, $1\le i\le n$ is less 
than $10 C_0 |\tilde{Q}_i|$.
If $\tilde{Q}'$ is another cube of the smallest scale in $\tilde{\Q}$, 
then the convex hull $h(\overline{Q})$ 
has distance at least $\K |\tilde{Q}_i|$ from $h(\overline{Q}')$
by \eqref{sparse-single-scale}. Hence every interval
of length $20C_0 |\tilde{Q}_i|$ contains at least one
interval of length $|\tilde{Q_i}|$
which does not intersect $h(\overline{Q}')$
for any $\tilde{Q}'$ at the smallest scale.

Now consider a cube $\tilde Q$ at the second smallest scale.
We may increase each 
interval $1000 \tilde Q_i$ by at most $1$ percent on either side, 
thus obtaining an interval $\overline{Q}_i$, so that for any cube 
$\tilde{Q}'$ at the smallest scale we either have 
$h(\overline{Q}')\subseteq \overline{Q}_i$
or $h(\overline{Q}')\cap \overline{Q}_i=\emptyset$.
All we need is that the endpoints of the enlarged interval
are not contained in the convex hull $h(\overline{Q}')$
for any $\tilde{Q}'$ at the smallest scale, which can
clearly be achieved from \eqref{sparse-scale} and the previous discussion.

Then we may define the convex hull $h(\overline{Q})$
as before, and the separation properties of these convex 
hulls discussed for the smallest scale also hold for cubes at 
the second smallest scale.

By successively passing to larger scales, we can increase
each interval $1000 \tilde{Q}_i$ by at most $1$ percent on either side
so that it contains either all of $h(\overline{Q}')$
for any $\tilde{Q}'$ at a smaller
scale or is disjoint from $h(\overline{Q}')$. 
Namely each interval of length $|\tilde{Q}_i|$ contains an interval 
of length comparable to the sidelength of cubes at the next smaller scale
which does not intersect $h(\overline{Q}')$ for any $\tilde{Q}'$
at the next smaller scale,
and so on passing to smaller and smaller scales, 
so we can find an endpoint for $\overline{Q}_i$ which is not
contained in $h(\overline{Q}')$ for any $\tilde{Q}'$ of any smaller scale.

We can now formulate the good dyadic properties that the intervals 
$\overline{Q}_i$ have:

\begin{lemma}\label{dyadic-lemma}

Let $\tilde Q,\tilde Q'\in \tilde\Q$.
If $\diam(\tilde Q)<\diam(\tilde{Q}')$, then $10\overline{Q}_i\cap {\overline{Q}'}_j\neq 0$
for some $1\le i,j\le n$ implies 
$10\overline{Q}_{i'}\subseteq {\overline{Q}'}_j$ for all
$1\le i'\le n$.

\end{lemma}

\begin{proof}
The proof is clear by construction.  
\end{proof}

Let $\Q$ denote the set of boxes
$$ \Q := \{ L(\tilde Q): \tilde Q \in \tilde \Q \}.$$
We let $j_Q$ denote the rational number such that $2^{\K j_Q} = |Q_n|$, and refer to 
$j_Q$ as the \emph{frequency parameter} of $Q$.  
By \eqref{sparse-scale} the set of $j_Q$ is a discrete set so that
any two elements have at least distance $1$. 
For purely notational convenience we shall assume that $j_Q$
is an integer for each $Q\in \Q$. This is only a special case of the general situation, in which one may assume the $j_Q$ belong to a fixed 
shift of the integer lattice,
but the proof is the same in the general situation with some additional notation.

For each $Q \in \Q$ 
define $\pi_{Q_i}$ to be the Fourier multiplier
$$ \widehat{\pi_{Q_i} f}(\xi) := m_{L^{-1}(Q),i}(v_i^{-1}\xi) \hat f(\xi).$$
The notation $\pi_{Q_i}$ is a bit sloppy: to be more precise one should write
$\pi_{Q_i,i}$ or $\pi_{Q,i}$.
However, the index $i$ will always either appear in the subscript of $\pi$
or otherwise be clear from the context.
Note that the symbol of $\pi_{Q_i}$ is a bump function adapted to $Q_i$.  

By Plancherel, we can then rewrite the desired estimate in physical space as
$$
|\sum_{Q \in \Q} \int \prod_{i=1}^n
\pi_{Q_i} f_i(x)\ dx|
\leq C \prod_{i=1}^n \|f_i\|_{p_i}\ \ \ .$$

This completes our frequency space decomposition.  
Observe that we have returned to a notation that does not explicitly
involve any quantities with tilde accent. We shall not need the tilde 
accent anymore to denote quantities under the transformation $L^{-1}$, 
and thus be free to use the tilde accent in other contexts.

\section{Discretization in the physical space}

For fixed $Q$, the projections $\pi_{Q_i}$ are Fourier multipliers supported on intervals of length ranging from $|Q_n|$ to $2^{M_1} |Q_n|$.  The 
Heisenberg uncertainty principle then suggests that one needs to consider several spatial scales, from the coarse scale of $|Q_n|^{-1}$ to the fine scale of $2^{-M_1} |Q_n|^{-1}$.  This multiplicity of scales causes much technical difficulty in \cite{thiele}, \cite{grafakos-li}.
A key difference and simplification 
in our approach is that we only decompose in the coarsest scale $|Q_n|^{-1}$, so that 
we will sometimes localize less than the uncertainty principle suggests. 
Unfortunately we have only been
able to make this simplification work in the $2 < p < \infty$ case, which is
why this paper is restricted to this range of exponents.

In frequency space we have used compactly supported
cutoff functions to decompose the multiplier. Consequently,
we will continue to use truncations which have compact support in 
frequency space and merely satisfy rapid decay estimates in physical space.
Becasue of this, we can use the standard dyadic grid to partition physical
space as opposed to the carefully constructed grid we use for frequency space.
An interval is called dyadic if it is
of the form $[2^kn,2^k(n+1)]$ with integers $k$ and $n$.

Following standard procedure, we shall index the space-frequency decomposition using tiles.

\begin{definition}\label{tile-def}
Let $1 \leq i \leq n$.  An $i$-tile is a rectangle
$P = I_P \times \omega_P$ with area $2^{M_i}$, $I_P$ a dyadic interval, and $\omega_P$ an interval in the mesh $\{ Q_i: Q \in \Q \}$.  A multi-tile is an $n$-tuple $\pv = (P_1, \ldots, P_n)$
such that each $P_i$ is an $i$-tile, the interval $I_{P_i} = I_\pv$ is independent of $i$, and such that the \emph{frequency box} $Q_\pv := \prod_{i=1}^n \omega_{P_i}$ of $\pv$ is an element of $\Q$.  The \emph{frequency parameter} $j_\pv$ of a multi-tile is defined by $j_\pv := j_{Q_\pv}$.  In particular, we have $|I_\pv| = 2^{-\K j_\pv}$.
A multi-tile $\pv'$ is called a translate of $\pv$, if $Q_\pv=Q_{\pv'}$.
If $P=I_P\times \omega_P$ is an $i$-tile and $\omega_P=Q_{i}$, we write
$\overline{\omega}_{P}$ for $\overline{Q}_i$.

\end{definition}

Since we assumed $|Q_n|$ to be an integral power of $2^\K$, 
we also have that $|I_P|$ is an integral power of $2^\K$ for all tiles $P$.

Note that $n$-tiles have area 1. In the Walsh analogue \cite{thiele2} they
correspond to one dimensional subspaces of $L^2(\R)$, and will thus be easier to handle than the other tiles. For instance, we shall be able to obtain good $L^p$ bounds on these tiles in addition to $L^2$ bounds thanks to Lemma \ref{bernstein}. The $i$-tiles have area larger than one, and due to
our goal to prove uniform estimates in the $M_i$, we do not have any
control over the area. In the Walsh model these tiles correspond to high dimensional
subspaces of $L^2(\R)$, and we can do little more than considering 
good $L^2$ estimates on them.

We proceed to define the cutoff operators in physical space.
Let $\eta$ denote a fixed positive function with total $L^1$- mass 1 and with Fourier transform supported in $[-2^{-2\K}, 2^{-2\K}]$, satisfying the pointwise estimates 
\begin{equation}\label{eta-bounds}
 C^{-1} (1 + |x|)^{-N^2} \leq \eta(x) \leq C (1 + |x|)^{-N^2}.
\end{equation}
Here $N$ is the previously chosen constant which controls decay in 
physical space.

Let $\eta_j$ denote the function $\eta_j(x) := 2^{-\K j}\eta(2^{-\K j}x)$.  For any subset $E$ of $\R$, denote by $\chi_E$ the characteristic function of $E$
and define the smoothed out characteristic function $\chi_{E,j}$ by
$$ \chi_{E,j} := \chi_E * \eta_{j}.$$
Note that 
\be{chi-split}
\chi_{\biguplus_{\alpha \in A} E_\alpha,j} = \sum_{\alpha \in A} \chi_{E_\alpha,j}.
\end{equation}
Here $\biguplus_{\alpha \in A} E_\alpha$ denotes the disjoint union of the $E_\alpha$; this is the same concept as $\bigcup_{\alpha \in A} E_\alpha$ but is only defined when the $E_\alpha$ are disjoint.
In particular, if $\R = \biguplus_{\alpha \in A} E_\alpha$ then
\be{i-decomp}
1 = \sum_{\alpha \in A} \chi_{E_\alpha,j}.
\end{equation}
Informally, $\chi_{E,j}$ is a frequency-localized approximation to $\chi_E$.  In fact we have the pointwise estimate
\be{e-diff}
|\chi_{E,j}(x) - \chi_E(x)| \le C (1 + 2^{\K j} \dist(x, \partial E))^{-N^2+1},
\end{equation}
where $\partial E$ is the topological boundary of $E$; this is easily verified by checking the cases $x \in E$ and $x \in \R \backslash E$ separately.

Now let $\Pv_0$ denote the space of all multi-tiles. 
From \eqref{i-decomp} we have the identity
$$
\sum_{Q \in \Q} \int \prod_{i=1}^n
\pi_{Q_i} f_i(x)\ dx
= \sum_{\pv \in \Pv_0} \int \chi_{I_\pv,j_\pv}(x)
\prod_{i=1}^n \pi_{\omega_{P_i}} f_i(x)\, dx.$$
It thus suffices to show that
\[
|\sum_{\pv \in \Pv_0}
\int \chi_{I_\pv,j_\pv} \prod_{i=1}^n \pi_{\omega_{P_i}} f_i| \leq C \prod_{i=1}^n \|f_i\|_{p_i}.
\]
Again, by a standard limiting argument it suffices to prove this
estimate where $\Pv_0$ has been replaced by the set of all tiles
$\pv\in \Pv_0$ such that $I_{\pv}\subseteq [-2^{\K k},2^{\K k}]$ for some fixed
large $k$, provided the constant in the estimate does not depend on $k$.
We shall fix such $k$ and denote this subset by $\Pv_1$. Thus $\Pv_1$ 
is a finite set of tiles (recall that the set $\tilde{\Q}$ 
of possible frequency boxes had been restricted to a finite set). 
We thus are aiming to show
\be{total-tile}
|\sum_{\pv \in \Pv_1}
\int \chi_{I_\pv,j_\pv} \prod_{i=1}^n \pi_{\omega_{P_i}} f_i| \leq C \prod_{i=1}^n \|f_i\|_{p_i}.
\end{equation}

\section{The geometry of tiles and trees}\label{tile-sec}

The following definition establishes an order relation on
the set of tiles.

\begin{definition}\label{tile-order-def}
Let $P$, $P'$ be $i$-tiles for some $1 \leq i \leq n$.
We say that $P \leq P'$ if $I_P \subseteq  I_{P'}$ and 
$\overline{\omega}_{P} \supseteq \overline{\omega}_{P'}$. 
If $\pv$ and $\pv'$ are two multi-tiles, we say that $\pv\le \pv'$
if there exists a $1\le i\le n$ such that $P_i\le P_i'$.
\end{definition}

Clearly the order of $i$-tiles is transitive. However, by the
good dyadic property of Lemma \ref{dyadic-lemma}, we also have

\begin{lemma}\label{transitive}
The order on multi-tiles is transitive, i.e., $\pv\le \pv'$ and $\pv'\le \pv''$
imply $\pv\le \pv''$.
\end{lemma}

\begin{proof}
Assume $\pv\neq \pv'$ and $\pv'\neq \pv''$. Then
$\pv_i\le \pv_i'$ and $\pv_j'\le \pv_j''$ imply 
$\overline{Q}_i \supset \overline{Q}_i'$
and $\overline{Q}_j' \supset \overline{Q}_j''$ with
strict containment by \eqref{sparse-single-scale}. By
Lemma \ref{dyadic-lemma} we conclude 
$\overline{Q}_i \supset \overline{Q}_j'$ and
$\overline{Q}_j' \supset \overline{Q}_i''$. This
gives $\pv\le \pv''$ as desired.
\end{proof}

As is standard in the theory, the main argument shall consist
in splitting the set $\Pv_1$ into smaller subsets, for which the
name {\it tree} has become standard. 
We call a dyadic interval $\K${\it - dyadic}, if it has length
$2^{\K j}$ for some integer $j$.
If $\xi\in \Gamma'$, $1\le i\le n$, and $I$ is a $\K$-dyadic interval,
we define
$$\omega_{i,\xi,I}:=[\xi_i-\frac 12 2^{M_i} |I|^{-1},\xi_i+\frac 12 2^{M_i}|I|^{-1}] $$
and 
$$\overline{\omega}_{\xi,I}:=[L^{-1}(\xi)_i-500 |I|^{-1},L^{-1}(\xi)_i+500 |I|^{-1}]\ \ \ . $$
Observe that the right hand side of the last display does not
depend on $i$, because $L^{-1}(\xi)$ is on the diagonal.
We shall not need any good dyadic properties for the intervals 
$\overline{\omega}_{\xi,I}$.

\begin{definition}\label{tree-def} Let $\xi\in \Gamma'$,
let $I$ be a $\K$- dyadic interval, and let $T$ be a set of multi-tiles.
The triple $(T,\xi,I)$ is called a tree, if $T$ is non-empty,
$I_\pv\subseteq I$ for each $\pv\in T$, 
and for each $\pv\in T$ there is a $1\le i\le n$ such that
$\overline{\omega}_{\xi,I} \subseteq \overline{\omega}_{P_i}$.

We shall write $\omega_{i,T}$ and 
$\overline{\omega}_{T}$ for $\omega_{i,\xi,I}$ and $\overline{\omega}_{\xi,I}$.
The data $(\xi,I)$ are called the \emph{top data} of the tree.

We will frequently call the set $T$ itself a tree, with the understanding 
that top data $(\xi_T,I_T)$ are associated to $T$.
If $T$ is a tree, we define the \emph{box set} to be the set 
$\Q_T := \{ Q_\pv: \pv \in T\}$.
For each $Q \in \Q_T$, we define the \emph{support} $E_{Q,T}$ of $Q$ to be the set
\be{eqt-def}
E_{Q,T} := \bigcup \{ I_\pv: \pv \in T, Q_\pv = Q\}.
\end{equation}
We say that two trees $T$, $T'$ are \emph{distinct} if $T$ and $T'$ 
have no tiles in common, that is $T \cap T' = \emptyset$.  
(We are reserving the term \emph{disjointness} for a stronger property, that the tiles in $T$ and $T'$ do not overlap). \end{definition}

If we say that one tree $T$ is contained in another tree $T'$, $T \subseteq T'$, this simply means that the multi-tiles in $T$ are also in $T'$, and does not imply any additional relation between the top data 
$(\xi_T,I_T)$ and $(\xi_{T'},I_{T'})$.

A main observation for trees is that the box set of a tree
can be parameterized by the frequency parameter. In this sense
trees make the connection to Littlewood Paley theory, in which
frequency boxes are essentially parameterized by their scale.

\begin{lemma}\label{boxset-lemma}
If $T$ is a tree, then for each $j$ there is at most one 
$Q\in \Q_T$ such that $j=j_Q$.
\end{lemma}

\begin{proof} Let $Q,Q'\in \Q_T$, then there are $i$ and $i'$
such that
$\overline{\omega}_{T} \subseteq \overline{Q}_i$
and
$\overline{\omega}_{T} \subseteq \overline{Q}_{i'}'$
Thus $\overline{Q}_i\cap \overline{Q}_{i'}'\neq \emptyset$.
This implies $Q=Q'$ or $j_Q\neq j_{Q'}$ by \eqref{sparse-single-scale}.
\end{proof}

If $T$ is a tree, we define the {\it scale set}
$$\J_T:=\{Pj_Q:Q\in \Q_T\}\ .$$
For
$j\in \J_T$ we define $Q^j$ to be the $Q\in \Q_T$ for which
$j_Q=j$.

The main tree selection algorithm will consist of 
a tree selection process as follows:

\begin{definition}
A tree selection process shall consist of choosing a tree $T_1$
from $\Pv_1$, then choosing a tree $T_2$ from $\Pv_1\setminus T_1$
and so on. I.e., at the $k$-th step we choose a tree $T_k$ from
$\Pv_1\setminus (T_1\cup\dots\cup T_{k-1})$. We shall refer to the
trees $T_k$ as the selected trees.
\end{definition}

For two reasons it will be necessary that at each step
these trees be as large as possible.

\begin{definition}\label{max-tree-def}
Consider a subset $\Pv$ of $\Pv_1$ and top data $(\xi,I)$ 
as in Definition \ref{tree-def}.
Then the maximal tree $T^*$ in $\Pv$ with top data $(\xi_{T^*},I_{T^*})=(\xi,I)$ 
is the set of all $\pv\in \Pv$ such that $I_\pv\subseteq I$ and
$\overline{\omega}_{T}\subseteq \overline{\omega}_{P_i}$.

A tree selection process is called greedy, if at the $k$-th step
the tree $T_k$ is maximal in $\Pv_1\setminus (T_1\cup\dots\cup T_{k-1})$.
\end{definition}

One of the reasons to run a greedy selection process is to
gain a nesting property described by the following lemma:

\begin{lemma}\label{simple-nesting}
Let $T$ be one of the selected trees of a greedy tree selection process.
If $Q,Q'\in \Q_T$ and $j_Q<j_{Q'}$, then $E_{Q,T}\supset E_{Q',T}$.
\end{lemma}

\begin{proof}
Suppose under the assumptions of the lemma we had
$E_{Q,T}\not\supset E_{Q',T}$, then there was a $\pv'\in T$
with $Q_{\pv'}=Q'$ such that $I_{\pv'}\not\subseteq E_{Q,T}$.
Pick any $\pv\in T$ such that $Q_{\pv}=Q$, and let $\pv''$ be
a translate of $\pv$ such that $I_{\pv'}\subseteq I_{\pv''}$.
Clearly $\pv''$ is an element of $\Pv_1$.
The multi-tile $\pv''$ cannot be in the tree $T$ because of
$I_{\pv'}\not\subseteq E_{Q,T}$, but its geometry would qualify it
to be in the tree, namely, $I_\pv''\subseteq I_T$ and 
$\overline{\omega}_{T}\subseteq \overline{\omega}_{P''_i}$
for some $i$.
Thus it must have been selected for a tree
$T''$ at a previous stage of the selection process. 
However, the geometry of $\pv'$ qualifies it to be in the same 
tree $T''$, namely, $I_{\pv'}\subseteq I_{\pv''}\subseteq I_{T''}$
and $\overline{\omega}_{T''}\subseteq \overline{\omega}_{P'_i}$
because $\overline{\omega}_{T''}\subseteq \overline{\omega}_{P''_j}$
for some $j$ and $\overline{\omega}_{P''_j}\subseteq \overline{\omega}_{P'_i}$
by Lemma \ref{dyadic-lemma} and the fact that 
$\overline{\omega}_{P''_i}$ and $\overline{\omega}_{P'_i}$ intersect.
This gives a contradiction to the maximality of $T''$.
\end{proof}

The nesting property of the above lemma implies the following bound on the
cardinality of the finite sets $\partial E_{Q^j,T}$.

\begin{lemma}\label{count}
Let $T$ be a selected tree of a greedy selection process.
Then 
$$ \sum_{j \in \J_T} 2^{-\K j} \# \partial E_{Q^j,T} \le C |I_T|\ \  .$$
\end{lemma}

This should be compared with the trivial bound 
$\# \partial E_{Q^j,T} \leq 2^{\K j}|I_T|$. 

\begin{proof}
Since $E_{Q^j,T}$ is a finite union of intervals, there are two types of points in $\partial E_{Q^j,T}$; those that are the left endpoint of connected components in $E_{Q^j,T}$, and those that are right endpoints.  Clearly it suffices to prove the bound for left endpoints only.

For each $j \in \J_T$ and each left endpoint $x$ of $E_{Q^j,T}$ consider the interval $(x - 2^{-\K j}, x - 2^{-\K j - 1})$.  
We claim that these intervals are pairwise disjoint. This claim implies the conclusion
of the lemma because these intervals are contained in $3I_T$.

To prove the claim, assume $(x - 2^{-\K j}, x - 2^{-\K j - 1})$ and
$(x' - 2^{-\K j'}, x - 2^{-\K j' - 1})$ have nonempty intersection. 
If $j=j'$, we necessarily have $x=x'$ because then
both $x$ and $x'$ are endpoints of dyadic intervals of length $2^{-\K j}$. 
Thus we can assume $j'<j$. Then $x$ is contained in the interior of the 
dyadic interval $I'$ of length $2^{-\K j'}$ with right endpoint $x'$.
However, the interior of $I'$ is disjoint from $E_{Q^{j'},T}$ and $x$ is contained
in $E_{Q^j,T}$, a contradiction to Lemma \ref{simple-nesting}.
This proves the claim.
\end{proof}

We will use Lemma \ref{count} to replace the spatial truncations in 
\eqref{total-tile} by certain variants of themselves when we sum
over the multi-tiles in a selected tree.

It will be convenient to replace the sets
$E_{Q^j,T}$ by variants $\tilde{E}_j$
which have better regularity properties:

\begin{definition}\label{hull-definition}
Let $T$ be a tree, and let $\I_T$ be the collection of all
maximal $\K$- dyadic intervals $I\subseteq I_T$ which have the property that
$3 I$ does not contain any of the intervals $I_\pv$ with $\pv\in T$.
For an integer $j$ with $2^{-\K j}\le |I_T|$ let $\tilde{E}_j$ be the union of all intervals $I$
in $\I_T$ such that $|I|< 2^{-\K j}$ (We emphasize that we have strict
inequality here, which will make the index $j$ most natural, as we can
see for example in the following lemma).
For an integer $j$ with $2^{-\K j}>|I_T|$ we define $\tilde{E}_j=\emptyset$.
\end{definition}

The sets $\tilde{E}_j$ obviously depend on the tree $T$, but we suppress
this dependence. The construction of the sets $E_j$ appears implicitly
in \cite{laceyt3}.

Clearly the intervals in $\I_T$ form a partition of $I_T$.
The nice regularity properties are stated in the following lemma:
\begin{lemma}\label{nicer-lemma}
Any two neighboring intervals in $\I_T$ differ by at most a factor
$2^\K$ in length. 

The set $\tilde{E}_j$ is a union of dyadic intervals of length 
$2^{-\K j}$ and contains $E_{Q^j,T}$ if $j\in \J_T$.
\end{lemma}

\begin{proof}
To prove the first statement, 
let $I$ and $I'$ be two neighboring intervals of $\I_T$ and
assume $I$ is larger. Let $I''$ be the dyadic interval which
contains $I'$ and has $2^{-\K}$ times the length of $I$.
We have to prove $I''\subseteq I'$. However, $3I''$ is contained in 
$3I$, und thus does not contain any interval $I_\pv$ with $\pv\in T$.
By maximality of $I'$ we have $I''\subseteq I'$. This proves the first statement
of the lemma.

To prove the second statement, let $I\subseteq I_T$ be any 
$\K$- dyadic interval of length $2^{-\K j}$ and observe the dichotomy that either
$I$ is contained in an interval of $\I_T$,
or $I$ is
partitioned into intervals of $\I_T$ which are strictly smaller than $I$.
The latter is necessarily the case if $I\cap \tilde{E}_j$ is
nonempty or if $I=I_\pv$ for some $\pv\in T$. This proves the second statement.

\end{proof}

\begin{lemma}\label{get-tile-lemma}
With the notation as above, if $I_0$ is a $\K$-dyadic interval 
of length $2^{-\K j_0}$  such that $3I_0\cap \tilde{E}_{j_0}\neq \emptyset$,
then there is a multi-tile $\pv\in T$ with $|I_\pv|\le |I_0|$ such that
$I_\pv\subseteq 10I_0$
\end{lemma}

\begin{proof}
There is a dyadic interval $I_1$ of same length as $I_0$ which
is contained in $3I_0$ and $\tilde{E}_{j_0}$. By definition
of $E_{j_0}$, $3I_1$ contains an $I_\pv$ for some multi-tile $\pv\in T$.
This together with \eqref{sparse-scale} proves the lemma.
\end{proof}

The sets $\tilde{E}_j$ are obviously nesting, hence we have
the following analogue of Lemma \ref{count}

\begin{lemma}\label{count-again}
Let $T$ be any tree. Then with the notation as above, 
$$ \sum_{j \in \J} 2^{-\K j} \# \partial \tilde {E}_j \le C |I_T|\ \  .$$
For each $j \in \J$, let $\Omega_j$ be the collection of connected components of $\tilde E_j$; thus $\Omega_j$ is a finite collection of intervals.  
For each $I \in \Omega_j$, let $x^l_I$ and $x^r_I$ denote the left and right endpoints of $I$, and let $I^l_j$ and $I^r_j$ denote the intervals
$$ I^l_j := (x^l_I - 2^{-\K (j+\m_i)-1}, x^l_I - 2^{-\K (j+\m_i) - 2})$$
$$ I^r_j := (x^r_I + 2^{-\K (j+\m_i) - 2}, x^r_I + 2^{-\K (j+\m_i)-1}).$$
Then the intervals $I^l_j$ are disjoint as $j$ varies in the integers
with $2^{-\K j}\le |I_T|$ and $I$ varies in $\Omega_j$, 

Moreover, if $I_j^l$ is an interval in 
the above collection, then the distance to the next
interval ${I'}_{j'}^l$ is at least $2^{-\K(j+2)}$.

Similar statements hold for the $I^r_j$.  
\end{lemma}

\begin{proof}
Most of the proof is exactly as in the proof of Lemma \ref{count},
the only new statement is the one on the distance between two
neighboring intervals.

Let $I_j^l$ be such an interval. It suffices to show
every different interval ${I'}_{j'}^l$ with $j'\ge j$ has distance at least
$2^{-\K(j+2)}$ from $I_j^l$. The case $j'\le j$ then follows
by symmetry.

The case $j=j'$ is easy, since the distance between $I^l_j$ and
${I}^l_{j'}$ is a multiple of $2^{-\K j}$.

Thus consider $j'>j$.
Let $\tilde{I}$ be the dyadic interval of length $2^{-\K j}$ 
which contains $I_j^l$ and let $\tilde{I}'$ be the dyadic
interval of length $2^{-\K(j+1)}$ whose left endpoint
is equal to the right endpoint of $I$. Clearly
$I$ is disjoint from $\tilde{E}_j$, whereas $I'$
is contained in $\tilde{E}_j$. The crucial observation
is that $I'$ cannot be contained in $\tilde{E}_{j+1}$, because
two neighboring intervals in $\I_T$ differ by at most a factor $2^{\K}$.

Thus the distance from $I_j^l$ to any point in $\tilde{E}_{j+1}$
is larger than $2^{-\K(j+1)}$, which proves the claim.

\end{proof}

While we will only consider trees selected by a greedy selection process,
we remark that the sets $\tilde{E}_j$ can be defined and satisfy
the above lemmata for arbitrary trees, which do not necessarily 
come from a greedy selection process and thus may not satisfy a 
nesting property as in Lemma \ref{simple-nesting}.

The regularity properties 
of the sets $\tilde{E}_j$  
will be used to construct phase plane
projection operators. It is worth mentioning that a similar notion
of regularity called convexity was used in \cite{thiele2} to construct phase
plane projections in the Walsh setting, although technically the
use of this type of convexity to construct phase plane projections 
is quite different in the current paper.

We introduce the notion of lacunarity in a tree:

\begin{definition}
An element  $\pv \in T$ is called $i$-lacunary if
$(\xi_{T})_i\not\in 2\omega_{P_i}$, and it is called
$i$-non-lacunary if $(\xi_{T})_i\in 2\omega_{P_i}$.
If $A$ is a subset of $\{1,\dots,n\}$, then $T_A$ is defined
to be the set of all elements in $T$ which are $i$-lacunary for all
$i\in A$ and $i$-non-lacunary for all $i\notin A$.
\end{definition}

Clearly, a tree can be written as the disjoint union
of subsets $T_A$ parameterized by all subsets $A\subseteq \{1,\dots,n\}$.
Each of the sets $T_A$ is either empty or again a tree with
the same top data as $T$.

We have
\begin{lemma}\label{heritage}
If $T$ is a tree and $A\subseteq \{1,\dots,n\}$, then for each $Q$
in $\Q_{T_A}$ we have $Q\in \Q_T$ and 
$E_{Q,T_A}=E_{Q,T}$.
\end{lemma}
\begin{proof}
This follows easily from the fact that lacunarity depends only
on the frequency intervals.
\end{proof}
As a consequence of this lemma, if $T$ is a selected tree as in
Lemma \ref{simple-nesting}, then whenever $T_A$ is non-empty
and thus a tree, it satisfies the analogues of Lemmata 
\ref{simple-nesting} and \ref{count}.

There is always at least one lacunary index:

\begin{lemma}\label{lacunary-lemma}
If $T$ is a tree, then $T_{\emptyset}$ is empty.
\end{lemma}
\begin{proof} If $\pv\in T_{\emptyset}$, 
then $L^{-1}(\xi_T)\in 2\overline{Q}_\pv$.
Since $\xi_T\in \Gamma'$ this contradicts
\eqref{whitney1}.
\end{proof}

We return to the second reason for running a greedy
selection process. It gives a certain strong disjointness property of 
multi-tiles of selected trees, as described by the following lemmata:

\begin{lemma}
Let $T$ and $T'$ be two selected trees of a greedy selection
process and assume that $T$ has been selected prior to $T'$.
Let $\pv\in T$ and $\pv'\in T'$ be such that 
we have 
\begin{equation}\label{non-empty-10q}
10\omega_{P_i} \cap 10\omega_{P'_i}\neq \emptyset
\end{equation}
and
\begin{equation}\label{size-comp}
|\omega_{P_i}| < |\omega_{P'_i}|
\end{equation}
for some $i$.
Then $I_{\pv'}\cap I_T=\emptyset$.
\end{lemma}

\begin{proof} 
Assume to get a contradiction that $I_{\pv'}\cap I_T\neq \emptyset$.
By size comparison we then necessarily have 
$I_{\pv'}\subseteq I_T$.

By \eqref{size-comp} and \eqref{sparse-scale} we have 
$100|\omega_{P_i}|< |\omega_{P'_i}|$. Hence 
\eqref{non-empty-10q} implies
$20\omega_{P_i} \subseteq 20\omega_{P'_i}\neq \emptyset$,
which implies
$\overline{\omega}_{P_i}\subseteq \overline{\omega}_{P'_i}$
(We made a point of not using Lemma \ref{dyadic-lemma} to conclude this).

Hence the geometry of the multi-tile $\pv'$ qualifies it to be
in the tree $T$. But it is not, because different selected trees have
no multi-tiles in common, so by maximality of $T$ it
must be an element of a tree that was selected prior to $T$.
This contradicts the assumed order of selection of $T$ and $T'$.
\end{proof}

The above lemma requires information about the order in which
trees have been selected. In our application this information will be
provided in the form described by the next lemma. If $\pv$ is a 
multi-tile, $s$ a
real number, $\xi\in \R^n$ and $1\le i\le n$, define the interval
\[\omega_{\pv,\xi,i,s}^+:=
[\xi_i+2^{-s}5000C_0 |\omega_{P_i}|, \xi_i+2^{2-s}5000C_0|\omega_{P_i}|]
\ \ . \]

\begin{lemma}\label{plus-separation-lemma}
Let $1\le i\le n$.
Let $\T$ be a subset of the set of selected trees of a greedy
selection process such that
$T,T'\in \T$ and $(\xi_{T})_i>(\xi_{T'})_i$ imply that 
$T$ has been selected
prior to $T'$. Let $s$ be some real number. 
Let $T,T'$ be two (not necessarily different) trees in $\T$
and assume $\pv\in T$ and $\pv'\in T'$ such that
\begin{equation}\label{non-empty-10omega}
10\omega_{P_i} \cap 10\omega_{P'_i}\neq \emptyset\ \ \ ,
\end{equation}
$$|\omega_{P_i}|< |\omega_{P'_i}|\ \ \ ,$$
\begin{equation}\label{non-empty-omega}
\omega_{\pv,\xi_T,i,s}^+ \cap \omega_{\pv',\xi_{T'},i,s}^+\neq \emptyset\ \ \ .
\end{equation}
Then $I_{\pv'}\cap I_T=\emptyset$ and 
in particular $I_{\pv'}\cap I_{\pv}=\emptyset$ and the trees $T$ and $T'$ are
indeed different.
\end{lemma}

\begin{proof} By the previous lemma we only need to prove that $T$ has
been selected prior to $T'$.
However, \eqref{non-empty-omega} together with
$$100|\omega_{P_i}|< |\omega_{P'_i}|$$
implies $(\xi_{T})_i>(\xi_{T'})_i$, which in turn implies that $T$ has been
selected prior to $T'$.
\end{proof}

Similarly, we can define
\[\omega_{\pv,\xi,i,s}^- :=
[\xi_i-2^{2-s}5000C_0|\omega_{P_i}|, \xi_i-2^{-s}5000C_0|\omega_{P_i}|]\]
and then have an analoguous lemma to \ref{plus-separation-lemma}
which we do not state explicitly.

We shall need another variant of this theme. For a tree $T$ define
\[\omega_{i,T,s}^+ :=
[(\xi_T)_i+2^{-s}10|\omega_{i,T}|, (\xi_T)_i+2^{2-s}10|\omega_{i,T}|]\ \ .\]

\begin{lemma}\label{tree-separation-lemma}
Let $1\le i\le n$.
Let $\T$ be a subset of the set of selected trees of a greedy
selection process such that
$T,T'\in \T$ and $(\xi_{T})_i>(\xi_{T'})_i$ imply that 
$T$ has been selected
prior to $T'$. 

Let $s$ be some real number. 
Let $T_0, T_1,T_2,T_3 \in \T$ and assume we have
for $j=1,2,3$ 
\begin{equation}\label{non-10-tree}
10\omega_{i,T_0} \cap 10\omega_{i,T_j}\neq \emptyset\ \ \ ,
\end{equation}
\begin{equation}\label{comp-size-tree}
|\omega_{i,T_0}|\le |\omega_{i,T_j}|\ \ \ ,
\end{equation}
\begin{equation}\label{non-tree}
\omega_{i,T_0,s}^+ \cap \omega_{i,T_j,s}^+\neq \emptyset\ \ \ .
\end{equation}
If $T_1,T_2,T_3$ are all different, then
$I_{T_1}\cap I_{T_2}\cap I_{T_3}=\emptyset$.

\end{lemma}

\begin{proof}

As before, we can use \eqref{non-tree} to conclude that 
if $|I_{T'}|>|I_{T''}|$, for $T',T'';\in \{T_1,T_2,T_3\}$, 
then $T'$ has been selected prior to $T''$. Thus we may assume
$|I_{T_1}|\ge |I_{T_2}|\ge |I_{T_3}|$ and $T_1,T_2,T_3$
have been selected in this order.

Assume to get a contradiction that  $I_{T_1}\cap I_{T_2}\cap I_{T_3}\neq \emptyset$, then we have by dyadicity
$I_{T_1}\supset I_{T_2} \supset I_{T_3}$.
Let $\pv$ be a multi-tile in $T_1 \cup T_2 \cup T_3$
for which $|\omega_{P_i}|$ is minimal. Let $T\in \{T_1,T_2,T_3\}$ 
be the tree so that $\pv\in T$.
If there
is another multi-tile $\pv'$ for which 
$|\omega_{P'_i}|=|\omega_{P_i}|$, then we conclude
by \eqref{sparse-single-scale} that 
$\omega_{P'_i}=\omega_{P_i}$. Hence $\pv$ and $\pv'$ are
in the same tree $T$.

Now let $\pv'\in T'$ for some $T'\in \{T_1,T_2,T_3\}$
with $|\omega_{P'_i}|>|\omega_{P_i}|$.
We aim to show $\pv'\in T_1$, which will prove the lemma,
because then one of the trees $T_2$ and $T_3$ has to be empty,
a contradiction.

Observe that by \eqref{non-10-tree} and \eqref{comp-size-tree}
we have
${\omega}_{i,T_0}\subseteq 20 {\omega}_{i,T_1}$.
By a similar argument for $T'$ we have
$20{\omega}_{i,T_1}\cap 20 {\omega}_{i,T'}\neq \emptyset$.
By assumption on the size of tree tops we have
${\omega}_{i,T_1}\subseteq 40 {\omega}_{i,T'}$.
If $|{\omega}_{i,T_1}|<|{\omega}_{i,T'}|$, then this implies
$\overline{\omega}_{T_1}\subseteq \overline{\omega}_{T'}$,
which in turn implies $\pv'\in T_1$.

We may thus assume $|{\omega}_{i,T_1}|=|{\omega}_{i,T'}|$.
Then we can merely conclude 
$\overline{\omega}_{T_1}\cap \overline{\omega}_{T'}\neq \emptyset$
and
$\overline{\omega}_{T_1}\subseteq 3\overline{\omega}_{T'}$.
By a similar argument with $T'$ in place of $T_1$ we have
$\overline{\omega}_{T'}\cap \overline{\omega}_{T}\neq \emptyset$
and 
$\overline{\omega}_{T'}\subseteq 3\overline{\omega}_{T}$.
But $\overline{\omega}_{T}\subseteq \overline{\omega}_{P_j}$
for some index $1\le j\le n$. Hence
$\overline{\omega}_{T'}\subseteq 3\overline{\omega}_{P_j}$.
However, for some possibly different index $l$ we
have $\overline{\omega}_{T'}\subseteq \overline{\omega}_{P'_l}$.
Hence $3\overline{\omega}_{P_j}$ and $\overline{\omega}_{P'_l}$
have non-empty intersection, and thus by Lemma \ref{dyadic-lemma}
we have $10\overline{\omega}_{P_j} \subseteq \overline{\omega}_{P'_l}$.
However, we have seen before
$\overline{\omega}_{T_1}\subseteq 3\overline{\omega}_{T'}$ and we
have $3\overline{\omega}_{T'}\subseteq 9\overline{\omega}_{P_j}$.
Hence $\overline{\omega}_{T_1}\subseteq \overline{\omega}_{P'_l}$,
which proves the claim.

\end{proof}

\section{A few remarks on smooth truncations and weights}

We pause to prove a few lemmata on weight functions
and smooth truncations of functions,
 which are best separated from the
main string of arguments so as to not slow the main argument down
by these technical lemmata.

Given an interval $I$, we define the approximate cutoff function 
$\tilde \chi_I$ as $$ \tilde \chi_I(x) := (1 + \frac{\dist(x,I)}{|I|})^{-1}.$$

We shall need the following lemma, which may be of independent interest.

\begin{lemma}\label{almost}
Let $\I$ be a finite collection of disjoint intervals, and suppose that for each $I$ one has an $L^2$ function $f_I$.  Then
$$ \| \sum_{I\in \I} |I|^{1/2} \tilde \chi_I f_I \|_2 
\leq C (\sum_{I\in \I} |I|)^{1/2} \sup_{I\in \I} \|f_I\|_2.$$
\end{lemma}

Note that this lemma would be automatic from Cauchy-Schwarz if $\sum_I \tilde \chi_I^2$ was uniformly bounded.  Unfortunately, this function is only in BMO, so one needs to work a little harder.  This lemma is a special case of a phase space Bessel inequality \eqref{t-size-alt} that we will need in the sequel.

\begin{proof}  We may assume that the $f_I$ are real and positive.
By estimating $\tilde \chi_I$ by a weighted sum of $\chi_{AI}$ for dyadic $A \geq 1$ (as before $AI$ has the same center as $I$ but $A$ times its length),
it suffices by the triangle inequality to show

\begin{lemma}\label{induct-almost}
Let $\I$ be as above, $f_I$ be real and positive, and let $A \geq 1$.  Then there exists an absolute constant $M > 0$ such that
\be{hyp}
\| \sum_{I\in \I} |I|^{1/2} \chi_{AI} f_I \|_2 \leq 
M A^{1/2} (\sum_{I\in \I} |I|)^{1/2}
\sup_{I\in \I} \|f_I\|_2.
\end{equation}
\end{lemma}

\begin{proof}
Fix $A$.  
Let $M_A$ be the best constant $M$ for which \eqref{hyp} obtains for all subsets of $\I$ (in place of $\I$ itself) and all $f_I$; our objective is 
to show that $M_A$ is bounded uniformly in $A$ and $\I$. 

By construction there exists a subset $\I'$ of $\I$ and $f_I$ such that
$$ \| \sum_{I\in \I'} |I|^{1/2} \chi_{AI} f_I \|_2 \gtrsim 
M_A A^{1/2} (\sum_{I\in \I'} |I|)^{1/2} \sup_{I\in \I'} \|f_I\|_2.$$
Squaring this we obtain
$$ \sum_I \sum_J \langle |I|^{1/2} \chi_{AI} f_I, |J|^{1/2} \chi_{AJ} f_J \rangle \gtrsim M_A^2 A (\sum_I |I|) \sup_I \|f_I\|_2^2.$$
One only has a contribution if $J \subseteq 5AI$ or $I \subseteq 5AJ$.  By symmetry and positivity we thus have
$$ \sum_I \sum_{J: J \subseteq 5AI} \langle |I|^{1/2} \chi_{AI} f_I, |J|^{1/2} \chi_{AJ} f_J \rangle \gtrsim M_A^2 A (\sum_I |I|) \sup_I \|f_I\|_2^2.$$
Moving the $J$ summation inside the inner product, and using Cauchy-Schwarz and the definition of $M_A$, we may estimate the left-hand side by
$$ \sum_I |I|^{1/2} \|f_I\|_2 M_A A^{1/2} (\sum_{J: J \subseteq 5AI} |J|)^{1/2} 
\sup_J \|f_J\|_2.$$
On the other hand, by disjointness we have $\sum_{J: J \subseteq 5AI} |J| \leq 5A|I|$.  Inserting this into the previous estimates we obtain $M_A \leq C 1$
as desired.
\end{proof}
\end{proof}

Lemma \ref{almost} will be used through the following variant:

\begin{lemma}\label{almost-useful}
 Let $I'$ be an interval and let $\I$ be a collection of
disjoint intervals such that $|I|\le |I'|$ for each $I\in \I$.
For each $I\in \I$ consider a positive $L^2$ function $f_I$. Then
$$\| \tilde{\chi}_{I'}
 \sum_{I\in \I} |I|^{1/2} \tilde \chi_I f_I 
\|_2 \leq C |I'|^{1/2}       \sup_I \|f_I\|_2.$$
\end{lemma}

\begin{proof}

We may write
\[\tilde{\chi}_{I'}\leq C 
\chi_{I'}+\sum_{k\ge 0}2^{-k}
\chi_{2^{k+1}I'\setminus 2^{k}I'}\]
By the triangle inequality it suffices to prove for each $k\ge 0$
$$\| 
 \sum_{I\in \I:I\cap 2^{k+1}I'\setminus 2^{k}I'\neq \emptyset} 
|I|^{1/2} \tilde \chi_I f_I 
\|_2 \leq C 2^{k/2}
|I'|^{1/2}       \sup_I \|f_I\|_2.$$
and an analoguous estimate for the term $\chi_{I'}$, which is clearly an
easy variant of the above.
Since $|I|\le |I'|$ for all $I\in \I$ and the intervals $I$ are pairwise disjoint, we have
$$ \sum_{I\in \I:I\cap 2^{k+1}I'\setminus 2^{k}I'\neq \emptyset} \leq C
 2^k|I'|\ \ .$$
The claim now follows by Lemma \ref{almost}.
\end{proof}

For any symbol $m$ on $\R$, we let $T_m$ be the associated Fourier multiplier.  (This is consistent with \eqref{lambda-t} if we lift $m$ from $\R$ to $\{ (\xi,-\xi): \xi \in \R\}$ in the obvious manner).

Let $w$ be a positive function on $\R$ and $r > 0$ be a number.  We say that $w$ is \emph{essentially constant at scale $r$} if one has
$$ (1 + \frac{|x-y|}{r})^{-100} \leq C \frac{w(x)}{w(y)}
\leq C (1 + \frac{|x-y|}{r})^{100}$$
for all $x, y \in \R$.  In particular, the weights $\tilde \chi_I^\alpha$ are essentially constant at any scale $|I|$ or less when $|\alpha| \leq 100$.

We shall need the following weighted version of Bernstein's inequality.

\begin{lemma}\label{bernstein}
Let $f$ be a function whose Fourier transform is supported on an interval $\omega$ of width $O(2^{\K j})$ for some integer $j$. Then we have
$$ \| w f \|_\infty \leq C 2^{\K j/2} \| w f \|_2$$
for all weights $w$ which are essentially constant at scale $2^{-\K j}$.
\end{lemma}

\begin{proof}
We can write $f = T_m f$ where $m$ is a suitable bump function adapted to $2\omega$.  From the decay of the kernel of $T_m$ we thus have the pointwise estimate
$$ |f(x)| = |T_m f(x)| \leq C 2^{\K j} \int \frac{f(y)}{(1 + 2^{\K j} |x-y|)^{N}}\ dy$$
and the claim easily follows.
\end{proof}

Let $T$ be a (possibly vector-valued) convolution operator and $r > 0$ be a number.  We say that $T$ is \emph{essentially local at scale $r$} if the convolution kernel $K(x)$ satisfies the bounds
\be{rapdec}
|K(x)| \leq C r^N |x|^{-N}
\end{equation}
for $|x| \gg r$.

\begin{lemma}\label{local}  Let $T$ be a convolution operator which is bounded on $L^2$ and which is essentially local at some scale $r > 0$.  Then one has
\be{local-est}
\| w Tf \|_2 \leq C \| wf \|_2
\end{equation}
for all weights $w$ which are essentially constant at scale $r$.
\end{lemma}

\begin{proof}
We can truncate $T$ so that the kernel is supported on the interval $\{ |x| \leq C r \}$; from \eqref{rapdec} it is easy to see that this does not affect the $L^2$ boundedness of $T$ or \eqref{local-est}.  The claim then follows by partitioning space into intervals of length $r$ and applying the $L^2$ boundedness hypothesis to each interval separately.
\end{proof}

\section{Phase space norms and the size of a tree}\label{phase-sec}

The general approach to proving an estimate such as
\eqref{total-tile} is to prove the estimate first in the easier
case when the summation goes only over a tree rather
than the whole set $\Pv_1$. The point being that this easier 
estimate is a matter of standard Littlewood-Paley theory 
without modulation 
invariance. The top of the tree fixes a frequency,
which after a modulation we can think of as being the zero frequency 
in standard Littlewood-Paley theory. 
In our situation this Littlewood-Paley estimate is Proposition \ref{paraproduct-cor}.

The second step is to organize the whole set into trees by a greedy
selection process, so as to sum the tree estimates. In this organization, 
the notion of size of a tree plays a crucial role. In this section we shall
introduce this notion.

\begin{definition}\label{size-def}
In the following definitions $1 \leq i \leq n$, and $f_i$ is an $L^2$ function.

If $P_i$ is an $i$-tile and $\xi_i \in \R$, we define the semi-norm $\|f_i\|_{P_i,\xi_i}$ by
\be{fips}
\| f_i \|_{P_i, \xi_i} := \sup_{m_{P_i}} \| \tilde \chi_{I_\pv}^{10} T_{m_{P_i}}(f_i) \|_2
\end{equation}
where $m_{P_i}$ ranges over all smooth functions supported on $10 \omega_{P_i}$
which satisfy the estimates
\be{mpi-est}
|\partial_\xi^k m_{P_i}(\xi)| \leq |\xi - \xi_i|^{-k} \frac{|\xi - \xi_i|}{|\omega_{P_i}|}  
\end{equation}
for all $\xi \in \R$ and $0 \leq k \leq N^2$.  In particular, $m_{P_i}$ vanishes at $\xi_i$.

If $T$ is a tree, we define the \emph{$i$-size} $\size_{i}(T)$ of $T$ by
\be{size-nonmax-def}
\size_{i}(T) := (\frac{1}{|I_T|} \sum_{\pv \in T} \| f_i\|_{P_i,(\xi_T)_i}^2)^{1/2}+ 
|I_T|^{-\frac 12} 
\sup_{m_{i,T}}\|\tilde{\chi}_{I_T}^{10}T_{m_{i,T}}(f_i)\|_2.
\end{equation}
where $m_{i,T}$ ranges over all smooth functions supported on 
$10 \omega_{i,T}$
which satisfy the estimates
\be{mpi-est-tree}
|\partial_\xi^k m_{i,T}(\xi)| \leq |\xi - (\xi_T)_i|^{-k} 
\frac{|\xi - (\xi_T)_i|}{|\omega_{i,T}|}  
\end{equation}
for all $\xi \in \R$ and $0 \leq k \leq N^2$.

If $\Pv$ is any collection of multi-tiles, we define the maximal size $\size^*_{i}(\Pv)$ of $\Pv$ to be
\be{size-max-def}
\size^*_{i}(\Pv) := \sup_{(T,\xi,I):T \subseteq \Pv} \size_{i}(T)
\end{equation}
where $(T,\xi,I)$ ranges over all trees with $T\subseteq \Pv$.
\end{definition}

We remark, that for $|\xi-\xi_i| \gg |\omega_{P_i}|$
we have the following estimates for $m_{P_i}$ which are stronger
than \eqref{mpi-est}:
\be{mpi-further-est}
|\partial_\xi^k m_{P_i}(\xi)| 
\leq C |\xi - \xi_i|^{-k} (\frac{|\xi - \xi_i|}{|\omega_{P_i}|})^{1-l}\ .  
\end{equation}
for all $0\le l,k < N^2/2$.
These estimates can be obtained from \eqref{mpi-est} and the support condition on
$m_{P_i}$ by integrating the $m_{P_i}$ over
its support $l$ times. Thus $m_{P_i}$ is forced to be rather
small if $\xi_i$ is far away from its support. This observation shall however 
only be of technical importance, since in our applications $\xi_i$
will always be within $C\omega_{P_i}$ for some moderate constant $C$.

Thus, heuristically, $\|f_i\|_{P_i,\xi_i}$ is the $L^2$ norm of $f_i$ when restricted to the portion of $P_i$ which is away from the frequency $\xi_i$.
The tree size  $\size_{i}(T)$ is heuristically 
the $L^2$ norm of $f_i$ when restricted to the region $\bigcup_{\pv \in T} P_i\cup (I_T\times \omega_{i,T})$ in phase space.  As a gross caricature, one has the very approximate relationship
$$ \size_i(T) ``\approx'' \osc_{I_T}(e^{-2\pi i (\xi_T)_i \cdot} f_i)$$
although this caricature does not fully capture the phase space localization to $T$ in \eqref{size-nonmax-def}.

Here we have written $\osc_I(f)$ 
for the $L^2$ mean oscillation given by
\be{osc-def}
\osc_I(f) := (\frac{1}{|I|} \int_I |f - \frac{1}{I} \int_I f|^2)^{1/2}.
\end{equation}

If $\Pv$ is a tree, the maximal size $\size^*_i(\Pv)$ 
is a strengthened version 
of the tree size $\size_i(\Pv)$.  
The relationship between the two is analoguous 
to the relationship between the BMO norm on an interval $I$ and the $L^2$ 
mean oscillation on that interval $I$.  This analogy is particularly 
accurate when $M_i = 0$ and $(\xi_T)_i = 0$.

The freedom to choose $I$ in \eqref{size-max-def} independently of the
set $T$ adds a useful ``Hardy-Littlewood maximal function''-component
to the size definition. For example this is used in the following lemma:

\begin{lemma}
We have
\be{unsigned-est}
\| \tilde \chi_{I_\pv}^{10} T_{m_{P_i}}(f_i) \|_2
\leq C \| f_i \|_{P_i, \xi_i}
 \leq C |I_\pv|^{1/2} \size^*_i(T) 
\end{equation}
for all indices $1 \leq i \leq n$, trees $T$, multi-tiles $\pv \in T$, frequencies $\xi_i \in \R$, and symbols $m_{P_i}$ supported on $10 \omega_{P_i}$ which obey \eqref{mpi-est}. 

Moreover,
\be{unsigned-est-sup}
\| \tilde \chi_{I}^{10} T_{m_{i, I}}(f_i) \|_2
 \leq C |I|^{1/2} \size^*_i(T) 
\end{equation}
for all indices $1 \leq i \leq n$, trees $T$, 
all $i$-non-lacunary multi-tiles $\pv \in T$, 
all $\K$-dyadic intervals $I$ with $I_\pv\subseteq 10 I$ and all
symbols $m_{i,I}$ supported on $5\omega_{i,\xi_T,I}$ and satisfying
$$
|\partial_\xi^k m_{i,I}(\xi)| \leq {|\omega_{i,\xi_T,I}|^{-k}}  
$$
for all $\xi \in \R$ and $0 \leq k \leq N^2$.  
\end{lemma}

\begin{proof}
We first consider \eqref{unsigned-est}.
The first inequality is just by definition.
By the remark just after Definition \ref{size-def}, 
we only have to prove the second inequality for $\xi_i\in  100\omega_{P_i}$. 
Then this inequality follows from \eqref{size-nonmax-def}, 
\eqref{size-max-def} since $T$ contains the singleton tree $\{\pv\}$ with
top data $\xi',I'$ where $\xi'$ is defined by 
$\xi_i=(\xi')_i$ and $I'$ is the $\K$-dyadic
interval of length $2^\K |I_\pv|$ which contains $I_\pv$. 
It is easily verified that these top data indeed turn
$\{\pv\}$ into a tree.

Now we consider \eqref{unsigned-est-sup}.
Observe first of all that $|I_\pv|\le |I|$
because both intervals are $\K$-dyadic. By translating $I$, we may
as well assume $I_{\pv}\subseteq I$. Namely, we have to
translate $I$ by at most ten times its length, so that
$\tilde{\chi}_I$ stays the same up to some bounded factor.

We consider the two cases $|I_\pv|=|I|$ and
$|I_\pv|<|I|$. Assume first $|I_\pv|=|I|$.

Observe that by $i$-non-lacunarity, $\omega_{P_i}$
is contained in $5\omega_{i,\xi_T,I}$.
Hence $5\omega_{i,\xi_T,I}$ is contained in
$9 \omega_{P_i}$.
Pick $\xi'\in \Gamma'$ so that $(\xi')_i$ is an endpoint
of $10\omega_{P_i}$. Then
$$\| \tilde \chi_{I_\pv}^{10} T_{m_{i, I}}(f_i) \|_2\le 
C \|f_i\|_{P_i,(\xi')_i} $$
by definition of the right hand side, because the
multiplier $m_{i,I}$ is supported in $10\omega_{P_i}$ and
satisfies \eqref{mpi-est}, possibly with a constant. 
Now the claim follows from \eqref{unsigned-est},
which proves \eqref{unsigned-est-sup} in the case $|I_\pv|=|I|$.

Now assume $|I_\pv|<|I|$. We consider again a singleton
tree $\{\pv\}$ with top data $\xi$, $I$ so that 
$\xi_i$ is an endpoint of $5\omega_{i,\xi_T,I}$.
Again, by $i$-non-lacunarity of $\pv$ we see that
these top data indeed turn $\{\pv\}$ into a tree.
Then the multiplier $m_{i,I}$ satisfies \eqref{mpi-est-tree}
with respect to this singleton tree as one can easily see, 
possibly with some constant, and \eqref{unsigned-est-sup} follows 
by definition of the tree size.

\end{proof}

Inequality \eqref{unsigned-est-sup} will be mainly applied in
connection with Lemma \ref{get-tile-lemma}.

Note that our definition of size is $L^2$-based (but normalized to have an $L^\infty$ scaling, see the proposition below). In principle one can define
$L^p$ based notions of size by using the $L^p$ norm instead of the $L^2$ norm
in \eqref{fips} and the properly adjusted normalizations thereafter.
In the $M_i = 0$ case the $L^2$-based notion of size is essentially equivalent to the $L^p$-based notion thanks to Lemma \ref{bernstein} and the fact that the $P_i$ then have area 1. Only if $n=i$
this will be guaranteed, which is why this case will have special treatment.

However when $M_i \gg 1$ the $L^2$ based notion of size is not equivalent
to an $L^p$ based notion.  In order to make the Bessel 
inequality \eqref{t-size-alt} work we need $L^2$-based sizes, and this is one of the reasons for the restriction $2 < p_i < \infty$ in Theorem \ref{main}.  Presumably one would need arguments such as those in \cite{thiele} to remove this restriction.

We conclude this section with the observation that the size is always controlled by the $L^\infty$ norm.

\begin{proposition}\label{initial}  For all $1 \leq i \leq n$ and arbitrary functions $f_i$, we have
\be{size-initial}
\size^*_i(\Pv_1) \leq C \|f_i\|_\infty.
\end{equation}
\end{proposition}

\begin{proof}
Fix $i$.  It suffices to show that
$$\size_i(T) \leq C \|f_i\|_\infty $$
for all trees $T$ in $\Pv_1$.

Fix $T$.  
We have to estimate both summands in \eqref{size-nonmax-def}.
The second summand is immediate since $T_{m_{i,T}}$ is a universally bounded operator
in $L^\infty$.

We consider the first summand in \eqref{size-nonmax-def}.
By frequency modulation invariance we may assume that $\xi_T = 0$.  From \eqref{size-nonmax-def}, \eqref{fips} it suffices to show that
$$
\sum_{\pv \in T} \| \tilde \chi_{I_\pv}^{10} T_{m_{P_i}}(f_i) \|_2^2 \leq C C |I_T| \|f_i\|_\infty^2
$$
whenever $m_{P_i}$ is supported on $10 \omega_{P_i}$ and obeys \eqref{mpi-est} with $\xi_i = 0$.

First suppose that $f_i$ vanishes on $3I_T$.  Then a simple computation  
following Lemma \ref{local} 
shows that
$$
\| \tilde \chi_{I_\pv}^{10} T_{m_{P_i}}(f_i) \|_2^2
\leq C |I_\pv| (\frac{|I_\pv|}{|I_T|})^2 \|f_i\|_\infty^2$$
for all $\pv \in \T$, and the claim follows by summing in $\pv$.  Thus we may assume that $f_i$ is supported in $3I_T$.  It then suffices to show that
$$
\sum_{\pv \in T} \| \tilde \chi_{I_\pv}^{10} T_{m_{P_i}}(f_i) \|_2^2 \leq C C \| f_i \|_2^2.
$$
By Khinchin's inequality, we may estimate the left-hand side by
$$
\| \sum_{\pv \in T} \epsilon_{\pv} \tilde \chi_{I_\pv}^{10} T_{m_{P_i}}(f_i) \|_2^2$$
for a suitable choice of signs $\epsilon_\pv \in \{+1, -1\}$.  But then the claim follows since the expression inside the norm is simply a pseudo-differential operator of order 0 in the symbol class $S^{1,-1}$ applied to $f_i$, with bounds uniform in $T$, $M_i$, $m_{P_i}$, and $\epsilon_\pv$ (see \cite{stein}).
\end{proof}

We remark that one can improve the bound from $1$ to $|E_i|/|E_k|$ for some other $1 \leq k \leq n$, if one is willing to remove from $E_k$ the exceptional set where the Hardy-Littlewood maximal  function of $E_i$ is $\gg |E_i|/|E_k|$, and also to remove from $\Pv_1$ all multi-tiles whose spatial interval is contained in this exceptional set (cf. \cite{mtt-1} Lemma 7.8; see also \cite{thiele}, \cite{laceyt2}, \cite{grafakos-li}).  Such improved estimates will not be needed 
here, because we restrict attention to the case $p_i>2$.

\section{Tree estimates}\label{tree-sec}

To follow the approach outlined at the beginning of Section \ref{phase-sec},
it is necessary to estimate the sum on the left hand side
of \eqref{total-tile}, where the summation
over $\Pv_1$ is replaced by a summation over a tree, 
by the various $i$-sizes of the tree. In doing so we may assume that the
tree has been selected by a greedy selection process, because all our
trees will be selected this way.
This tree estimate will be Proposition \ref{tree-est-prop} 
below.

\begin{proposition}\label{tree-est-prop}  Let $T$ be a selected tree of
a greedy selection process, and let $f_1, \ldots, f_n$ be test 
functions on $\R$ which satisfy the normalization
\be{infty-bound}
\| f_i \|_\infty \le 1
\end{equation}
for all $1 \leq i \leq n$.  Then we have
\be{tree-est}
|\sum_{\pv \in T}
\int \chi_{I_\pv,j_\pv} \prod_{i=1}^n \pi_{\omega_{P_i}} f_i| \leq 
C_{(\theta_i)}
|I_T| \prod_{i=1}^n \size^*_i(T)^{\theta_i}
\end{equation}
whenever $\theta_n = 1$, and $0 < \theta_1, \ldots, \theta_{n-1} < 1$ are such that $\sum_{i=1}^{n-1} \theta_i < 2$, and the implicitly used constant
$N$ is sufficiently large depending on $(\theta_i)$.
\end{proposition}

In spirit this estimate is a version of Proposition \ref{paraproduct-cor} 
with additional localization in physical space. It would be very tedious
to adapt the proof of Proposition \ref{paraproduct-cor} to the localized
setting, we shall therefore follow a different approach.
The idea is to replace the $f_i$
by functions (phase plane projections) which are very close to 
$f_i$ near the phase plane region of the tree, but 
essentially
vanish outside this region.
Then we shall apply Proposition \ref{paraproduct-cor} to these
phase plane projections directly and recover the result
of Proposition \ref{tree-est-prop} from there: the tree sizes
on the right hand side of \eqref{tree-est} will essentially be
the suitably normalized $L^p$-norms of the phase plane projections.

The parameters $\theta_i$ will be chosen depending on the exponents $p_i$
in Theorem \ref{main}. In the following we shall not write explicitly
$(\theta_i)$-dependence of our constants, in the same way as we
do not write $p_i$-dependence explicitly.

In analogy to \cite{laceyt1}, \cite{laceyt2}, \cite{mtt-1} one might expect a bound of $|I_T| \prod_{i=1}^n \size^*_i(T)$ on the right-hand side of \eqref{tree-est}.  This bound is achievable in the non-degenerate case $\m_i = O(1)$; however, in general only the $n^{th}$ size $\size_n^*(T)$ can be recovered
with a full power $\theta_n=1$.  
This gain in the case $i=n$ will be
crucial in the rest of the paper. One should probably be able to obtain the endpoint $\sum_{i=1}^{n-1} \theta_i = 2$ of this result, 
but we shall not attempt to do so here.

We will prove Proposition \ref{tree-est-prop} in this and the next section.
We remark that knowledge of the proof of \eqref{tree-est} will
not be needed in the later sections.

\begin{proof} 
Fix $T$, $\theta_i$, $f_i$.
We shall exploit scale invariance to reduce to the case $|I_T| = 1$.  
By a frequency modulation leaving $\Gamma$ and $\Gamma'$
invariant, we may as well assume $\xi_T = 0$.  

By the triangle inequality it suffices to prove estimate \eqref{tree-est}
where the summation goes over the subset $T_A$ instead of $T$ for
any subset $A\subseteq \{1,\dots,n\}$. Fix such $A$ and let
$B$ denote the set of non-lacunary indices, i.e., 
$\{1, \ldots, n\} = A \uplus B$.
We may assume $T_A$
is non-empty and thus a tree with top data $(\xi_T,I_T)$.
By Lemma \ref{lacunary-lemma} we have that $A$ is non-empty. 
Changing notation we may as well
assume that all multi-tiles in $T$ have lacunarity type $A$ and thus $T=T_A$.
Recall that $T_A$ also satisfies the crucial nesting property of Lemma
\ref{simple-nesting}.

Let $\J$ denote the set $\J := \{j_Q: Q \in \Q_T\}$; its smallest
element is greater or equal $0$.
Recall that by Lemma 
\ref{boxset-lemma}
the map $Q \mapsto j_Q$ is one-to-one. 

From \eqref{chi-split} we may estimate the left-hand side of \eqref{tree-est} as
\be{tree-1}
\sum_{j \in \J} |\int \tilde \chi_j \prod_{i=1}^n \tilde \pi_j f_i|
\end{equation}
where 
the spatial cutoff $\tilde \chi_{j}$ and Fourier cutoff $\tilde \pi_j$ are defined as
\be{tchij-def}
\tilde \chi_{j} := \chi_{E_{Q^j,T},j} = \sum_{\pv \in T: Q_\pv = Q^j} \chi_{I_\pv,j}
\end{equation}
and
$$ \tilde \pi_j := \pi_{Q^j_i}$$
for $j \in \J$.  In particular, $\tilde \chi_{j}$ has Fourier support in the region $\{\xi:  |\xi| \le 2^{\K j-\K} \}$.  
We remark that our notation is sloppy here, the operator $\tilde{\pi}_j$ also
depends on the parameter $i$. We will always write $\tilde{\pi}_j$
in combination with a function $f_i$, and the omitted index is
always the one of the function $f_i$.

Let $2 < p_i < \infty$ be chosen so that $p_i < 2/\theta_i$ for $
1\le i\le n-1$ and $\sum_{i=1}^n 1/p_i = 1$; this can be done by the hypotheses on $\theta$ stated in Proposition \ref{tree-est-prop}.  
We observe the simple bounds

\begin{lemma}\label{triv-bounds}
For all $1 \leq i \leq n$, $j \in \J$, and intervals $I$ of length $|I| = 2^{-\K j}$ we have
$$ \| \tilde \chi_j^{1/2n} \tilde \pi_j f_i \|_{L^{p_i}(I)} \leq C
|I|^{1/p_i} \size^*_i(T)^{\theta_i}.$$
\end{lemma}

\begin{proof}
By interpolation and Proposition \ref{initial} it suffices to prove the bounds
$$ \| \tilde \chi_j^{1/2n} \tilde \pi_j f_i \|_{L^2(I)} \leq C
|I|^{1/2} \size^*_i(T),$$
$$ \| \tilde \chi_j^{1/2n} \tilde \pi_j f_i \|_{L^{\infty}(I)} \leq C,$$
and (in the $i=n$ case only)
$$ \| \tilde \chi_j^{1/2n} \tilde \pi_j f_i \|_{L^{\infty}(I)} \leq C \size^*_i(T).$$
The second estimate is immediate from the boundedness of the $f_i$, while the third follows from the first and Lemma \ref{bernstein} since $|v_n|=1$ and thus $\m_n = 0$.  
Thus it suffices to prove the first inequality.  

Fix $i$, $j$, $I$.  From \eqref{tchij-def} we see that there exists $\pv \in \T$ with $|I_\pv|=|I|$ 
such that we have the pointwise estimate
$$ \tilde \chi_j(x)^{1/2n} \leq C \tilde{\chi}_{I_\pv}^{10} $$ 
on $I$.  It thus suffices to show that
$$  \| \tilde \chi_{I_\pv}^{10} \tilde \pi_j f_i \|_{L^2} \leq C |I|^{1/2} \size^*_i(T).$$
This however was observed in \eqref{unsigned-est}.
\end{proof}

From this Lemma and H\"older we have
$$ \int_I |\tilde \chi_j \prod_{i=1}^n \tilde \pi_j f_i|
\leq C \|\tilde \chi_j\|_{L^\infty(I)}^{1/2} |I| \prod_{i=1}^n \size^*_i(T)^{\theta_i}$$
for all $j \in \J$ and $|I| = 2^{-\K j}$.  Summing in $I$, we see that we may bound each summand in \eqref{tree-1} by the right-hand side of \eqref{tree-est}.  Thus the main difficulty is to obtain summability in $j$.

In the non-degenerate case $\m_i = O(1)$ this summability is obtained by noting there must be at least two lacunary indices\footnote{This is due to the fact that the cubes $10 Q^j$ must intersect $\hyperplane$ in order to have a non-zero contribution to \eqref{tree-1}.} in $A$, estimating those indices in $L^2$ in space (and $l^2$ in $j$), and taking all other indices in $L^\infty$.  See e.g. \cite{laceyt1}, \cite{laceyt2}, \cite{mtt-1} and the discussion in the introduction of \cite{mtt-2}. This however is not feasible in the general case for two reasons.  Firstly, there might only be one lacunary index; and secondly, one can only get good $L^\infty$ bounds for the $i$ index when $\m_i = O(1)$.   
Thus we will have to invoke Proposition \ref{paraproduct-cor} as outlined 
before.

We shall need to replace $\chi_j$ in \eqref{tree-1} by a product of cutoff
functions. If $\chi_j$ was a characteristic function, we could simply
write it as a power of itself. However, it is an approximate characteristic
function. From \eqref{e-diff} we have:
\be{tchi-est}
|\tilde \chi_j(x) - \chi_{E_{Q^j,T}}|\le  C (1 + 2^{\K j} \dist(x, \partial E_{Q^j,T}))^{-N^2+1} .
\end{equation}

The right hand side can be handled by Lemma \ref{count}.

\begin{lemma}\label{triv-est}
Let $1\le k\le n$, then
$$
\sum_{j \in \J} |\int (\tilde \chi_j - \tilde \chi_j^{k})
\prod_{i=1}^n \tilde \pi_j f_i| \le
C |I_T|
\prod_{i=1}^n \size^*_i(T)^{\theta_i}.
$$
\end{lemma}

\begin{proof}
By the triangle inequality, we can estimate this sum by
$$
\sum_{j \in \J} \sum_{I: |I| = 2^{-\K j}} \int_I  
|\tilde \chi_j - \tilde \chi_j^{k}| \prod_{i=1}^n |\tilde \pi_j f_i|.$$
By H\"older and Lemma \ref{triv-bounds}, this is bounded by
$$
\sum_{j \in \J} \sum_{I: |I| = 2^{-\K j}} |I| 
\|(\tilde \chi_j - \tilde \chi_j^{k})\tilde{\chi}_j^{-\frac 12}\|_{L^\infty(I)}
\prod_{i=1}^n \size^*_i(T)^{\theta_i}.$$
On the other hand, from \eqref{tchi-est} and the lower bound in \eqref{eta-bounds} we have
$$
|I| \|(\tilde \chi_j - \tilde \chi_j^{k})\tilde{\chi}_j^{-\frac 12}\|_{L^\infty(I)}
\leq C \int_I (1 + 2^{\K j} \dist(x, \partial E_{Q^j,T}))^{-N}.
$$
Summing this in $I$ and evaluating the integral, we estimate the previous by
$$ C \sum_{j \in \J} 2^{-\K j} \# \partial E_{Q^j,T}$$
and the claim follows from Lemma \ref{count}.
\end{proof}

Instead of estimating \eqref{tree-1}, it suffices by the above Lemma
to to show that
\be{tree-2}
\sum_{j \in \J} |\int 
\prod_{i=1}^n \tilde{\chi}_j \tilde{\pi}_j f_i|
\le C
|I_T| \prod_{i=1}^n \size^*_i(T)^{\theta_i}.
\end{equation}

To do this we shall now introduce the crucial tree projection operators.  

\begin{proposition}\label{proj}
Let the notation and assumptions be as above.  Then for each $1 \leq i \leq n$ there exists a function $\Pi_i(f_i)$ with the following properties.
\begin{itemize}
\item (Control by size)
We have
\be{fire-1}
\| \Pi_i(f_i) \|_{p_i} \leq C |I_T|^{ 1/{p_i}}\size^*_i(T)^{\theta_i}.
\end{equation}
\item ($\Pi_i(f_i)$ approximates $f_i$ on $T$: lacunary case)
For all $i \in A$ and $j \in \J$, we have
\be{fire-2}
\tilde \chi_j \tilde \pi_j f_i = S_{j+ m_i} \Pi_i(f_i)
\end{equation}
where $S_{j+ \m_i}$ is a suitable Littlewood-Paley  projection to the frequency region $\{\xi: \frac 1{10} 2^{\K (j+\m_i)}\le  |\xi| \le 5000C_0 2^{\K (j+\m_i)}\}$.

\item ($\Pi_i(f_i)$ approximates $f_i$ on $T$: non-lacunary case) 
For all $i \in B$, $j_0 \in \J$, and all intervals $I_0$ of length $2^{-\K j_0}$, we have the bounds
\be{error-1}
\| \tilde \chi_{j_0}^{1/2n} \tilde \pi_{j_0} (f_i - \Pi_i(f_i)) \|_{L^{p_i}(I_0)}
\leq C \size^*_i(T)^{\theta_i} |I_0|^{1/p_i-1} \int \tilde \chi_{I_0}^2 \mu_{j_0}
\end{equation}
where $\mu_j$ is the function
\be{muj-def}
\mu_j(x) := \sum_{0 \leq j' } 2^{-|j'-j|/100} \sum_{y \in \partial \tilde E_{j'}} (1 + 2^{\K j'} |x-y|)^{-100}
\end{equation}
and the sets $\tilde{E}_j$ have been 
defined in Definition \ref{hull-definition}.
Also, for all $1 \leq i \leq n$, $j \in \J$, and intervals $I$ of length $|I| = 2^{-\K j}$ we have
\be{fire-loc}
\| \tilde \chi_j^{1/2n} \tilde \pi_j \Pi_i(f_i) \|_{L^{p_i}(I)} \leq C
|I|^{1/p_i} \size^*_i(T)^{\theta_i}.
\end{equation}
\end{itemize}
\end{proposition}

One can construct the $\Pi_i$ to be linear operators, but we shall not use this.

Roughly speaking, $\Pi_i(f_i)$ is the projection of $f_i$ to the region $\bigcup_{\pv \in T} P_i$ of phase space.  Although such a description can easily be made rigorous in a Walsh model, and is not too difficult in a lacunary Fourier model, it is substantially more delicate in the Fourier setting in the non-lacunary case.

Note that, the function $\mu_j$ can be controlled by Lemma \ref{count-again}.

We prove this rather technical Proposition in Section \ref{tech-sec}.  For now, we see why this proposition, combined with Proposition \ref{paraproduct-cor}, gives \eqref{tree-est}.

To exploit the fact that $\Pi_i(f_i)$ approximates $f_i$ on $T$, we use \eqref{fire-2} and the triangle inequality to estimate \eqref{tree-2} by the main term
\be{tree-main}
\sum_{j \in \J} |\int  
(\prod_{i \in A} \tilde \chi_j \tilde \pi_j f_i)
(\prod_{i \in B} \tilde \chi_j \tilde \pi_j \Pi_i(f_i))|
\end{equation}
plus $\# B$ error terms of the form
\be{tree-error}
\sum_{j \in \J} |\int \tilde \chi_j \tilde{\pi}_j (f_{i_0} - \Pi_{i_0}(f_{i_0}))
 (\prod_{i \in A \hbox{ or } i < i_0} \tilde{\chi}_j\tilde \pi_j f_i) 
\prod_{i \in B: i > i_0} \tilde{\chi}_j\tilde \pi_j \Pi_i(f_i)|
\end{equation}
where $i_0 \in B$.

Let us first control the contribution of a term \eqref{tree-error}.  Fix $i_0$. By the triangle inequality we may estimate \eqref{tree-error} by 
$$
\sum_{j \in \J} \sum_{I: |I| = 2^{-\K j}}
\int_I 
|\tilde \chi_j \tilde{\pi}_j (f_{i_0} - \Pi_{i_0}(f_{i_0}))|
 (\prod_{i \in A \hbox{ or } i < i_0} |\tilde{\chi}_j \tilde \pi_j f_i|)
\prod_{i \in B: i > i_0} |\tilde{\chi}_j
\tilde \pi_j \Pi_i(f_i)|.
$$

By H\"older, we may estimate the previous by
$$
\sum_{j \in \J} \sum_{I: |I| = 2^{-\K j}}
\| \tilde \chi_j \tilde {\pi}_j (f_{i_0} - \Pi_{i_0}(f_{i_0})) \|_{L^{p_{i_0}}(I)}\times$$
\be{fly}
\times \prod_{i \in A \hbox{ or } i < i_0} \| \tilde{\chi}_j\tilde \pi_j f_i \|_{L^{p_i}(I)}
\prod_{i \in B: i > i_0} \| \tilde{\chi}_j\tilde \pi_j \Pi_i(f_i) \|_{L^{p_i}(I)}.
\end{equation}
The first factor we estimate using \eqref{error-1}.  The second group of factors we estimate using Lemma \ref{triv-bounds}.  
For the last group we use \eqref{fire-loc}.  Combining all these estimates and using \eqref{scaling}, we see that we may estimate \eqref{fly} by
$$
C
\sum_{j \in \J} \sum_{I: |I| = 2^{-\K j}}
(\prod_{i=1}^n \size^*_i(T)^{\theta_i}) \int \tilde \chi_I^2 \mu_j$$
which of course sums to
$$
C
(\prod_{i=1}^n \size^*_i(T)^{\theta_i}) \sum_j \int \mu_j.$$
Expanding out $\mu_j$, we may estimate this by
$$
C (\prod_{i=1}^n \size^*_i(T)^{\theta_i}) \sum_j 
\sum_{y \in \partial \tilde E_j} \int (1 + 2^{\K j} |x-y|)^{-100}\ dx;$$
computing the integral, we thus obtain
$$
C (\prod_{i=1}^n \size^*_i(T)^{\theta_i}) \sum_j 
2^{-\K j} \#\partial \tilde E_j$$
and the claim follows from Lemma \ref{count-again}.

It remains to estimate \eqref{tree-main}.  By repeating the proof of Lemma \ref{triv-est} (but using \eqref{fire-loc} in place of Lemma \ref{triv-bounds} when $i \in B$) it suffices to estimate
$$
\sum_{j \in \J} |\int  
(\prod_{i \in A} \tilde \chi_j \tilde \pi_j f_i)
(\prod_{i \in B} \tilde \pi_j \Pi_i(f_i))|
$$
which we rewrite using \eqref{fire-2} as
$$
\sum_{j \in \J} |\int  
(\prod_{i \in A} S_{j + \m_i} \Pi_i(f_i))
(\prod_{i \in B} \tilde \pi_j \Pi_i(f_i))|.
$$
By Proposition \ref{paraproduct-cor}, we may estimate this by
$$C 
\prod_{i = 1}^n \|\Pi_i(f_i)\|_{p_i},$$
and the claim follows from \eqref{fire-1}. Observe that in applying Proposition
\ref{paraproduct-cor} we have used $A\neq \emptyset$, which had been a consequence of Lemma
\ref{lacunary-lemma}.

This concludes the proof of Proposition \ref{tree-est-prop}, except for the proof of Proposition \ref{proj} which shall be done in the next section.
\end{proof}

\section{Proof of Proposition \ref{proj}}\label{tech-sec}

We now prove Proposition \ref{proj}.  The proof of this Proposition is rather involved, and the techniques used here are not needed elsewhere in the paper; readers who are interested in the general shape of the proof of Theorem \ref{main} may wish to skip this section on a first reading.  

We begin with the lacunary case $i \in A$, which is substantially easier.

Fix $i \in A$.  We define $\Pi_i(f_i)$ as
\be{F-def}
\Pi_i(f_i) := \sum_{j \in \J} \tilde \chi_j \tilde \pi_j f_i
\end{equation}
in this case.  From the Fourier support of $\tilde \chi_j$ and $\tilde \pi_j$ we see that \eqref{fire-2} is obeyed.  It remains to show \eqref{fire-1}.

By interpolation and Proposition \ref{initial} it suffices to prove the estimates
\be{2-est}
\| \Pi_i(f_i)  \|_2 \leq C |I_T|^{1/2}\size^*_i(T)
\end{equation}
\be{q-est}
\| \Pi_i(f_i)   \|_{\BMO} \leq C
\end{equation}
with the additional estimate
\be{bmo-est}
\| \Pi_i(f_i)   \|_{\BMO} \leq C \size^*_i(T)
\end{equation}
when $i=n$. Here we read $\BMO$ as $\K$- dyadic $\BMO$ (see the proof below).

We first prove \eqref{2-est}.  By orthogonality and \eqref{tchij-def} we have
$$ \| \Pi_i(f_i) \|_2^2 = \sum_{j \in \J} \| \sum_{\pv \in T: Q_\pv = Q^j} \chi_{I_\pv,j} \tilde \pi_j f_i \|_2^2.$$
Since we have the pointwise bound
$$ \sum_{\pv \in T: Q_\pv = Q^j} \chi_{I_\pv,j} \leq C$$
we can bound the previous using Cauchy-Schwarz by
$$ C \sum_{j \in \J} \sum_{\pv \in T: Q_\pv = Q^j} \| |\chi_{I_\pv,j}|^{1/2} \tilde \pi_j f_i \|_2^2.$$
But by \eqref{unsigned-est} this is bounded by
$$ C \sum_{j \in \J} \sum_{\pv \in T: Q_\pv = Q^j} \| f_i \|_{P_i,0}^2$$
and the claim follows from \eqref{size-nonmax-def}.

The claim \eqref{q-est} follows immediately from the observation that the linear operator $\Pi_i$ is a pseudo-differential operator of order 0 in the symbol class $S^{1,-1}$, and therefore maps $L^\infty$ to BMO (see e.g. \cite{stein}), so it remains to show \eqref{bmo-est}.

We need to show that
\be{osci}
\osc_I(\Pi_n(f_n)) \leq C \size^*_n(T)
\end{equation}
for all $\K$-dyadic intervals $I$, which is what we mean by
$\K$- dyadic $\BMO$.

Fix $I$. We expand the left-hand side of \eqref{osci} using \eqref{tchij-def} as
$$
\osc_I(\sum_{\pv \in \T} \chi_{I_\pv,j_\pv} \pi_{\omega_{P_n}} f_n).
$$

First consider the contribution of the coarse scales, when $|I_\pv| > |I|$.  By the sub-linearity of $\osc_I()$, it suffices to show that
\be{big-tiles}
\sum_{\pv \in \T: |I_\pv| > |I|} \osc_I(\chi_{I_\pv,j_\pv} \pi_{\omega_{P_n}} f_n) \leq C \size^*_n(T).
\end{equation}
We use the easily verified Poincare inequality
$$ \osc_I(f) \leq C (|I|  \int_I |\nabla f|^2)^{1/2}$$
The Fourier multipliers $|I_\pv| \nabla \pi_{\omega_{P_n}}$ and $\pi_{\omega_{P_n}}$
have symbols adapted to $10 \omega_{P_n}$ which vanish at the origin,
whereas $\chi_{I_\pv,j_\pv}$ and $I_\pv \nabla \chi_{I_\pv,j_\pv}$ are dominated by
$\tilde{\chi}_{I_\pv}^{200}$.
By the previous and \eqref{unsigned-est}, we thus have
$$ \osc_I(\chi_{I_\pv,j_\pv} \pi_{\omega_{P_n}} f_n) \leq
C  
\frac{|I|}{|I_\pv|} (1 + \frac{\dist(I,I_\pv)}{|I_\pv|})^{-100}
\size^*_n(T).$$
Summing over all $\pv$ such that $|I_\pv| > |I|$ we see that this contribution is acceptable.

It remains to consider the contribution of the fine scales, i.e. those $\pv$ in the set 
$$\T_I := \{ \pv \in T: |I_\pv| \leq |I| \}.$$
For this contribution we will shall use the estimate
$$ \osc_I(f) \leq C \frac{1}{|I|^{1/2}} \| f \chi_{I,j_I} \|_2,$$
where $2^{\K j_I} := |I|$.  It thus suffices to show that
\be{oro}
\|\sum_{\pv \in \T_I} \chi_{I,j_I} \chi_{I_\pv,j_\pv}  \pi_{\omega_{P_n}} f_n\|_2 \leq C |I|^{1/2} \size^*_n(T).
\end{equation}
For fixed $j_\pv$, the expression inside the norm has Fourier support in the region $|\xi| \sim 2^{\K j_\pv}$ (lacunary supports).  
By orthogonality we may thus estimate the left-hand side of \eqref{oro} by
$$
C (\sum_{j \in \J: j \geq j_I} \|\sum_{\pv \in \T_I: j_\pv = j} \chi_{I,j_I} \chi_{I_\pv,j_\pv}  \pi_{\omega_{P_n}} f_n\|_2^2)^{1/2}.
$$
For fixed $j_\pv = j$, the $\chi_{I_\pv,j_\pv}(x)$ are uniformly summable in $I_\pv$ and $x$. 
By Cauchy-Schwarz we may thus estimate the previous by
$$
C (\sum_{j \in \J:j \geq j_I} \sum_{\pv \in \T_I: j_\pv = j} \|\chi_{I,j_I} \chi_{I_\pv,j_\pv}^{1/2} \pi_{\omega_{P_n}} f_n\|_2^2)^{1/2}.
$$
On the other hand, from \eqref{unsigned-est} we see that
$$
\|\chi_{I,j_I} \chi_{I_\pv,j_\pv}^{1/2} \tilde \pi_{\omega_{P_n}} f_n\|_2
\leq C (1 + \frac{\dist(I, I_\pv)}{|I|})^{-100} \| f_n \|_{P_n,0}.$$
Thus we may estimate the left-hand side of \eqref{oro} by
$$
C (\sum_{\pv \in \T_I} 
(1 + \frac{\dist(I, I_\pv)}{|I|})^{-200} \| f_n \|_{P_n,0}^2)^{1/2}.
$$
We now break up this as a sum over sub-trees of $\T$.
Let $\T^*_I$ denote those multi-tiles in $\T_I$ for which the interval $I_\pv$ is maximal.  Clearly the $I_\pv$ are disjoint as $\pv$ varies over $\T^*_I$.  We can thus rewrite the previous as
$$
C (\sum_{\pv' \in \T^*_I} \sum_{\pv \in \T: I_\pv \subseteq  I_{\pv'}}
(1 + \frac{\dist(I, I_\pv)}{|I|})^{-200} \| f_n \|_{P_n,0}^2)^{1/2}.
$$
Since $|I_{\pv'}| \leq |I|$, we can estimate this by
$$
C (\sum_{\pv' \in \T^*_I} 
(1 + \frac{\dist(I, I_{\pv'})}{|I|})^{-200} 
\sum_{\pv \in \T: I_\pv \subseteq  I_{\pv'}}
\| f_n \|_{P_n,0}^2)^{1/2}.
$$
For each $\pv' \in \T^*_I$, the collection $\{ \pv \in \T: I_\pv \subseteq  I_{\pv'}\}$ is a subtree of $\T$ with top data $(\xi, I_{\pv'})$.  
By \eqref{size-max-def} we can thus estimate the previous by
$$
C (\sum_{\pv' \in \T^*_I} 
(1 + \frac{\dist(I, I_{\pv'})}{|I|})^{-200} |I_{\pv'}| \size^*_n(T)^2)^{1/2}.$$
Since the $I_\pv$ are disjoint and of size $\leq |I|$, we can bound this by
$$ C \size^*_n(T) (\int (1 + \frac{\dist(I,x)}{|I|})^{-200}\ dx)^{1/2}$$
and the claim follows.

It remains to handle the more difficult non-lacunary case.  Fix $i \in B$.

For each real number $j$, let $T_j$ be a Fourier multiplier (defined, say, by
dilations of a fixed multiplier) whose symbol is adapted to the frequency 
region 
$\{ |\xi_j| \le 2^{2+\K j}\}$, whose symbol equals 1 for $\{ |\xi_j| \le 2^{1+\K j} \}$, and let $S_j$ be the associated Littlewood-Paley  projections $S_j := T_j - T_{j-1}$. We may assume that
the kernels of $T_j$ and $S_j$ are real and even.

In \eqref{F-def}, one used lacunary Fourier multipliers followed by smooth spatial cutoffs to construct $\Pi_i(f_i)$.  In the non-lacunary case these smooth spatial cutoffs are not desirable as they interfere\footnote{Of course, one could try to control this error with commutator estimates, however one does not seem able to recover the crucial $2^{-|j-j_0|/100}$ decay in \eqref{muj-def} by this approach.} with the ability to decompose non-lacunary multipliers as a telescoping series of lacunary ones.  To avoid this difficulty we are forced to use rough spatial cutoffs instead.

A first guess as to the construction of $\Pi_i(f_i)$ would be
$$ \tilde \Pi_i(f_i) := \chi_{\tilde E_{0}} T_{j(x) + \m_i} f_i,$$
where for each $x \in \tilde E_{0}$ we define the integer-valued function $j(x)$ by
$$j(x) := \max \{ 0 \leq j : x \in \tilde E_{j} \}$$

One can expand $\tilde \Pi_i(f_i)$ as a telescoping series:
\be{telescope}
\tilde \Pi_i(f_i) = \chi_{\tilde E_0} T_{\m_i} f_i + \sum_{1 \leq j} \chi_{\tilde E_j} S_{j+\m_i} f_i.
\end{equation}
This proposed projection turns out to obey \eqref{fire-1}, but does not obey \eqref{error-1} due to the poor frequency localization properties of the characteristic functions $\chi_{\tilde E_j}$ in \eqref{telescope}.  Specifically, the cutoffs destroy the vanishing moments of the $S_{j+\m_i} f_i$, and this will cause a difficulty when trying to sum in $j$ because the projection $\tilde \pi_{j_0}$ is non-lacunary.  

To get around this problem we shall modify each term of $\tilde \Pi_i(f_i)$ (except for the first term $\chi_{\tilde{E}_0} T_{\m_i} f_i$) to have a zero mean.  In order that these modifications do not collide with each other, we shall place them in disjoint intervals. 

We recall the intervals $I_j^r$ and $I_j^l$ and the collections $\Omega_j$
introduced in Lemma \ref{count-again}.
Let $\phi^l_{I,j}$ and $\phi^r_{I,j}$ be bump functions adapted to $I^l_j$, $I^r_j$ with total mass 
$$ \int \phi^l_{I,j} = \int \phi^r_{I,j} = 2^{-\K (j + \m_i)}.$$

For each $j \ge 1$ and $I \in \Omega_j$, decompose $\chi_I$ as $\chi_I = H^l_I + H^r_I$, where $H^l_I(x) := H(x - x^l_I)$ and $H^r_I(x) := -H(x - x^r_I)$ are shifted Heaviside functions.

For each $j \ge 1$ and $I \in \Omega_j$, define the quantities $c^l_{I,j}$ and $c^r_{I,j}$ by
\be{cl-def}
c^l_{I,j} := 2^{\K (j+\m_i)} \int H^l_I S_{j+\m_i} f_i\ \ ,
\end{equation}
$$ c^r_{I,j} := 2^{\K (j+\m_i)} \int H^r_I S_{j+\m_i} f_i\ \ .$$

A basic estimate on $c^l_{I,j}$ is

\begin{lemma}\label{clij-est}
Let $j\ge 0$ and $I \in \Omega_j$.  
Then we have the estimate
\be{big-bound}
|c^l_{I,j}| \leq C 2^{\K (j+\m_i)} 
\int \frac{S_{j+\m_i} f_i(x)\ dx}{(1 + 2^{\K(j+\m_i)} |x-x^l_I|)^{100}}\ \ .
\end{equation}
In particular, we have
\be{clb-2}
|c^l_{I,j}| 
\leq C 2^{\K \m_i/2} \size^*_i(T)
\end{equation}
and
\be{clb}
|c^l_{I,j}| \leq C.
\end{equation}
\end{lemma}

\begin{proof}

From \eqref{cl-def} we may write 
$$ c^l_{I,j} = 2^{\K (j+\m_i)} \int (\tilde{S}_{j+\m_i} 
H^l_I) S_{j+\m_i} f_i\ .$$
Here $\tilde{S}_{j+\m_i}$ is a Littlewood Paley projection whose
multiplier is supported in $5\omega_{i,T}$, is constant $1$ for
$\xi\in 4\omega_{i,T}\setminus 2\K^{-1}\omega_{i,T}$, and vanishes on
$\K^{-1}\omega_{i,T}$.
The claim \eqref{big-bound} then follows by using repeated integration by 
parts to obtain pointwise bounds on $\tilde{S}_{j+\m_i} H^l_I$.

There exists a $\K$-dyadic interval $I'$ of length
$2^{-\K j}$ whose left endpoint coincides with $x^l_I$, because
$\tilde{E}_j$ and hence $I$ is a union of dyadic intervals of length
$2^{-\K  j}$ by Lemma \ref{nicer-lemma}. 
By Lemma \ref{get-tile-lemma} there is a tile $\pv\in T$
with $I_\pv\subseteq 10I'$.

The bounds \eqref{clb-2}, \eqref{clb} then follow from \eqref{big-bound},
\eqref{unsigned-est-sup}, \eqref{infty-bound}.
\end{proof}

We can now define the corrected projections $\Pi_i(f_i)$ as
\be{pif-def}
 \Pi_i(f_i) := \tilde \Pi(f_i) - \sum_{1 \leq j} \sum_{I \in \Omega_j} (c^l_{I,j} \phi^l_{I,j} + c^r_{I,j} \phi^r_{I,j}).
\end{equation}

We first verify \eqref{fire-1}.  We first observe that
\be{pif-bound}
\| \Pi_i(f_i) \|_\infty \leq C.
\end{equation}
Indeed, from \eqref{infty-bound} we see that $T_{j_0+\m_i} f_i$ and $\tilde \Pi_i(f_i)$ are bounded.  The remaining terms of \eqref{pif-bound} then follow from \eqref{clb} and the disjointness of the $\phi^l_{I,j}$, and similarly for $c^r_{I,j}$ and $\phi^r_{I,j}$.  

In light of \eqref{pif-bound} it suffices by interpolation and Proposition \ref{initial} to prove \eqref{2-est}, together with the bound
\be{gen-infty}
\| \Pi(f_i) \|_\infty \leq C 2^{\K \m_i/2} \size^*_i(T)\ \ .
\end{equation}
In fact, to prove \eqref{fire-1} we only need the bound \eqref{gen-infty} 
for $i=n$, but the general case will be useful later.

We now prove \eqref{gen-infty}.  From \eqref{clb-2} the contribution of the $c^l_{I,j}$ is acceptable.  Similarly for $c^r_{I,j}$.  Thus it remains to show that
$$ \| \tilde \Pi(f_i) \|_\infty \leq C 2^{\K \m_i/2} \size^*_i(T),$$
or in other words that
$$ |T_{j(x)+\m_i} f_i(x)| \leq C 2^{\K \m_i/2} \size^*_i(T)$$
for all $x \in \tilde E_0$.

Fix $x$, and define $j := j(x)$.  From the definition of $j(x)$ there exists 
a $\K$-dyadic interval $I'\subseteq \tilde{E}_j$ of length $2^{-\K j}$ which
contains $x$. 
It thus suffices to show that
$$ \| T_{j+\m_i} f_i \|_{L^\infty(I')} \leq C 2^{\K \m_i/2} \size^*_i(T).$$
By Lemma \ref{get-tile-lemma} there is a multi-tile $\pv \in T$ with 
$I_{\pv}\subseteq 10I'$. Hence the desired estimate
follows from \eqref{unsigned-est-sup} and Lemma \ref{bernstein}.  

This completes the proof of \eqref{gen-infty}.

It remains to prove \eqref{2-est}.  From \eqref{pif-def}, the triangle inequality, and the disjointness and size of the $\phi^l_{I,j}$ (and of the $\phi^r_{I,j}$) it suffices to show that
\be{t-2-est}
\| \tilde \Pi(f_i) \|_2 \leq C |I_T|^{1/2}\size^*_i(T)
\end{equation}
and
\be{c-sum}
(\sum_{0 \leq j}\sum_{I \in \Omega_j} 
|c^l_{I,j}|^2 2^{-\K(j + \m_i)} )^{1/2} \leq C |I_T|^{1/2}\size^*_i(T)
\end{equation}
together with a similar bound for the $c^r_{I,j}$.

The bound \eqref{c-sum} follows from \eqref{clb-2} and Lemma \ref{count-again}.

To prove \eqref{t-2-est} we expand
$$ \tilde \Pi(f_i) = \sum_{0 \leq j } \chi_{\tilde E_j \backslash \tilde E_{j+1}} T_{j + \m_i} f_i
= \sum_{0 \leq j } \sum_{I \cap \tilde E_j \backslash \tilde E_{j+1}\neq \emptyset: |I| = 2^{-\K j}} \chi_{I\cap \tilde E_j \backslash E_{j+1}} 
T_{j + \m_i} f_i.$$

As $j$ and $I$ vary in the above sum, the sets 
$I\cap \tilde E_j\backslash \tilde E_{j+1}$ are pairwise disjoint, hence it satisfies to show
$$ 
\sum_{0 \leq j } \sum_{I \cap \tilde E_j \backslash \tilde E_{j+1}\neq \emptyset: |I| = 2^{-\K j}} 
\| T_{j + \m_i} f_i \|_{L^2(I)}^2 \leq C |I_T| \size^*_i(T)^2$$

For each $j$ and $I$ in this sum there is a dyadic 
interval $I'\subseteq I$ of length $2^{-\K} |I|$ which is contained in
$\tilde E_j \backslash \tilde E_{j+1}$. 
This follows from Lemma \ref{nicer-lemma}.
As $j$ and $I$ vary, these intervals $I'$ are pairwise disjoint. Hence
it suffices to show that
$$ \| T_{j + \m_i} f_i \|_{L^2(I)} \leq C |I|^{1/2} \size^*_i(T)$$
for all $j$, $I$ in the above sum.  But for such $j$, $I$ we can find a multi-tile 
$\pv \in T$ with $I_\pv \subseteq 10I$ by Lemma \ref{get-tile-lemma}.  
The claim then follows from \eqref{unsigned-est-sup}.  
This proves \eqref{t-2-est} and thus \eqref{2-est}.
The proof of \eqref{fire-1} is now complete.

The estimate \eqref{fire-loc} will follow from \eqref{error-1}, Lemma \ref{triv-bounds}, the triangle inequality, and the fact that the $\mu_j$ are uniformly bounded.  Thus it only remains to verify \eqref{error-1}.  

Fix $j_0\ge 0$ and $I_0$ such that $|I_0| = 2^{-\K j_0}$.  From the frequency support of $\tilde \pi_{j_0}$ we may replace $f_i - \Pi_i(f_i)$ with $T_{j_0+\m_i} f_i - \Pi_i(f_i)$.  By interpolation and Proposition \ref{initial} it suffices to show the bounds
\be{mess-1}
\| \tilde \chi_{j_0}^{1/2n} \tilde \pi_{j_0} (T_{j_0+\m_i} f_i - \Pi_i(f_i)) \|_{L^\infty(I_0)}
\leq C \frac{1}{|I_0|} \int \tilde \chi_{I_0}^2 \mu_{j_0}
\end{equation}
and
\be{mess-2}
\| \tilde \chi_{j_0}^{1/2n} \tilde \pi_{j_0} (T_{j_0+\m_i} f_i - \Pi_i(f_i)) \|_{L^2(I_0)}
\leq C \size^*_i(T) |I_0|^{-1/2} \int \tilde \chi_{I_0}^2 \mu_{j_0}
\end{equation}
with the additional bound
\be{mess-3}
\| \tilde \chi_{j_0}^{1/2n} \tilde \pi_{j_0} (T_{j_0+\m_i} f_i - \Pi_i(f_i)) \|_{L^\infty(I_0)}
\leq C \size^*_i(T) \frac{1}{|I_0|} \int \tilde \chi_{I_0}^2 \mu_{j_0}
\end{equation}
when $i=n$.

We first show \eqref{mess-1}.  For this estimate the only bound we use on 
$f_i$ is \eqref{infty-bound}.

First suppose that $I_0$ is outside $\tilde E_{j_0}$.  Then the claim follows from \eqref{infty-bound}, \eqref{pif-bound}, the decay of $\tilde \chi_{j_0}$ and the estimate
\be{muj-lower}
(1 + 2^{\K j_0} \dist(I_0, \partial \tilde E_{j_0}))^{-N} \leq C \frac{1}{|I_0|} \int \tilde \chi_{I_0}^2 \mu_{j_0}.
\end{equation}

It remains to consider the case when $I_0$ is inside $\tilde E_{j_0}$.  In this case the cutoff $(\tilde \chi_{j_0})^{1/2n}$ is useless and will be discarded.  We now decompose
\begin{align}
T_{j_0+\m_i} f_i - \Pi_i(f_i) 
=& \chi_{\R \backslash \tilde E_{j_0}} T_{j_0+\m_i} f_i 
\label{d1}\\
&- \chi_{\R \backslash \tilde E_{j_0}} \tilde \Pi_i(f_i)
\label{d2}\\
&+ \chi_{\R \backslash \tilde E_{j_0}} 
\sum_{1 \leq j \leq j_0} \sum_{I \in \Omega_j} c^l_{I,j} \phi^l_{I,j}
\label{d3}\\
&+ \chi_{\R \backslash \tilde E_{j_0}} \sum_{1 \leq j \leq j_0} \sum_{I \in \Omega_j} c^r_{I,j} \phi^r_{I,j}
\label{d4}\\
&- \sum_{j_0 < j } \sum_{I \in \Omega_j} 
H^l_I S_{j+\m_i} f_i - c^l_{I,j} \phi^l_{I,j}\label{d5}\\
&- \sum_{j_0 < j } \sum_{I \in \Omega_j} 
H^r_I S_{j+\m_i} f_i - c^r_{I,j} \phi^r_{I,j}.\label{d6}
\end{align}

From \eqref{infty-bound} and the separation between $I_0$ and 
$\R \backslash \tilde{E}_{j_0}$ and the decay of the kernel of $\tilde \pi_{j_0}$ we can bound the contribution of \eqref{d1} by \eqref{muj-lower} as desired.  The terms \eqref{d2}, \eqref{d3}, \eqref{d4} are similar, although \eqref{d3}, \eqref{d4} use \eqref{clb} instead of \eqref{infty-bound}.

Now consider \eqref{d5}.  The idea is to interact the smoothing properties of $\tilde \pi_{j_0}$ with the moment property 
\be{mean}
\int H^l_I S_{j+\m_i} f_i - c^l_{I,j} \phi^l_{I,j} = 0
\end{equation}
coming from the construction of the $c^l_{I,j}$.

By the triangle inequality and the definition of $\mu_j$ it suffices to show that
\begin{equation}\label{infty-inside}
\| \tilde \pi_{j_0} (H^l_I S_{j+\m_i} f_i - c^l_{I,j} \phi^l_{I,j}) \|_{L^\infty(I_0)}
\leq C 2^{-(j-j_0)/100} 2^{-\K (j-j_0)} (1+\frac {\dist(I_0,x^l_I)}{|I_0|})^{-30}
\end{equation}
for all $j > j_0$ and $I \in \Omega_j$.

Fix $j, I$.  
We first assume $x^l_I \in 3I_0$,
in which case we may replace the right-hand side by $C 2^{-(j-j_0)/100} 2^{-\K(j-j_0)}$. 
Let $K(x)$ denote the convolution kernel of $\tilde \pi_{j_0}$.  
By \eqref{mean} 
it suffices to show that
$$ |\int (K(x_0 - x) - K(x_0 - x^l_I))
(H^l_I S_{j+\m_i} f_i(x) - c^l_{I,j} \phi^l_{I,j}(x))\ dx|
\leq C 2^{-|j-j_0|/100} 2^{-\K(j-j_0)}$$
for all $x_0 \in I_0$.

Fix $x_0\in I_0$.
The contribution of $c^l_{I,j} \phi^l_{I,j}$ can be controlled using 
\eqref{clb}, the fact that the support of
 $\phi^l_{I,j}$ is within a distance of $C2^{-\K(j+\m_i)}$ from $x_I^l$,
and the bound
\begin{equation}\label{infty-kernel-bound}
|K(x_0 - x) - K(x_0 - x^l_I)| \leq C 2^{2\K (j_0 + \m_i)} |x - x^l_I|.
\end{equation}

To deal with the $H^l_I S_{j + \m_i} f_i$ term we rewrite it as
$$ |\int (K(x_0 - x) - K(x_0 - x^l_I))
([H^l_I, S_{j+\m_i}] f_i(x))\ dx|$$
with the commutator 
$[H^l_I, S_{j+\m_i}]$
of $H^l_I$ and $S_{j+\m_i}$. Here we have used
that $\pi_{j_0}S_{j+\m_i}=0$ and the kernel of $S_{j+\m_i}$ has mean zero.
However, we have the easily verified estimate
$$|[H^l_I, S_{j+\m_i}]f_i(x)|\le  C (1+2^{\K (j+\m_i)}(x-x^l_I))^{-200} 
\|f_i\|_\infty\ \ \ .$$
The desired bound for the $H^l_I S_{j + \m_i} f_i$ term now follows
again from \eqref{infty-kernel-bound}.

It remains to consider the case $x^l_I \not \in 3I_0$.
This is done similarly to the previous case, however using 
instead of \eqref{infty-kernel-bound} the kernel bound
$$|K(x_0 - x) - K(x_0 - x^l_I)| \leq C 2^{2\K (j_0 + \m_i)} 
(1+2^{\K j}|x_0-x^l_I|)^{-200}|x - x^l_I|$$
for $2|x-x_I^l|<|x-x_0|$.

The treatment of \eqref{d6} is similar.  This concludes the proof of \eqref{mess-1} in all cases.

The inequality \eqref{mess-3} follows from \eqref{mess-2} and Lemma \ref{bernstein}.  Thus it only remains to show \eqref{mess-2}.

Now we show \eqref{mess-2}.  This will be a reprise of the proof of \eqref{mess-1}, except that we shall rely on \eqref{unsigned-est-sup} instead of \eqref{infty-bound}.

We turn to the details. 
We shall only concern ourselves with the case when $5I_0$ intersects $E_{j_0}$.
The case when $5I_0$ is disjoint from $E_{j_0}$ is done
with similar arguments as those below; one loses some powers of the separation between $I_0$ and $E_{j_0}$ whenever one applies \eqref{unsigned-est}, but this is more than compensated for by the decay of $\tilde \chi_{j_0}$, which also gives the additional factor of \eqref{muj-lower}.  We omit the details. 

Discarding the cutoff $\tilde \chi_{j_0}$, we reduce to
$$\| \tilde \pi_{j_0} (T_{j_0+\m_i} f_i - \Pi_i(f_i)) \|_{L^2(I_0)}
\leq C \size^*_i(T) |I_0|^{1/2}\int \tilde{\chi}^2_{I_0} \mu_{j_0}.
$$

Since $I_0$ is near $E_{j_0}$, we can find a 
multi-tile $\pv \in T$ such that $I_\pv \subseteq 10I_0$
by Lemma \ref{get-tile-lemma}.
From \eqref{unsigned-est-sup} we have
\be{suzerain}
\|\tilde \chi_{I_0}^{10} T_{j_0 + \m_i} f_i \|_2 \leq C |I_0|^{1/2} \size^*_i(T).
\end{equation}

Decompose $T_{j_0+\m_i} f_i - \Pi_i(f_i)$ as $\eqref{d1} - \eqref{d2} + \eqref{d3} + \eqref{d4} - \eqref{d5} - \eqref{d6}$ as before.

First consider the contributions \eqref{d1}, \eqref{d2}, \eqref{d3}, \eqref{d4} which come from outside $\tilde E_{j_0}$.  Let us first examine the non-local portion of these contributions, or more precisely
$$
\| \tilde \pi_{j_0} (\chi_{\R \backslash 3I_0} (\eqref{d1} - \eqref{d2} + \eqref{d3} + \eqref{d4})) \|_{L^2(I_0)}.$$
From the rapid decay of the kernel of $\tilde \pi_{j_0}$ we may estimate this by
$$ 2^{-100 \K \m_i} \| \tilde \chi_{I_0}^{100} (\eqref{d1} - \eqref{d2} + \eqref{d3} + \eqref{d4})) \|_{2}.$$
But each of these terms is acceptable thanks to \eqref{gen-infty} (or more precisely, the arguments used to prove \eqref{gen-infty} but applied to \eqref{d1}, \eqref{d2}, \eqref{d3}, \eqref{d4} rather than $\Pi_i(f_i)$).

It remains to estimate the local contribution 
$$
\| \tilde \pi_{j_0} (\chi_{3I_0} (\eqref{d1} - \eqref{d2} + \eqref{d3} + \eqref{d4})) \|_{L^2(I_0)}.$$

Of course, these contributions are non-zero only when $3I_0$ is not contained in $\tilde E_{j_0}$, so that $I_0$ is near the boundary of $\tilde E_{j_0}$.  We can then discard the projection $\tilde \pi_{j_0}$ and the $\mu$- factor
on the right hand side and reduce to showing that
$$ \| \eqref{d1} - \eqref{d2} + \eqref{d3} + \eqref{d4} \|_{L^2(3I_0)} 
\leq C |I_0|^{1/2} \size^*_i(T).$$

The contribution of \eqref{d1} is acceptable by \eqref{suzerain}.

Now consider \eqref{d2}.  The set $(\tilde E_0 \backslash \tilde E_{j_0}) \cap 3I_0$ is the union of at most three intervals $I_1$ of length $2^{-\K j_0}$.  On each of these intervals $I_1$, the function $j(x)$ is equal to $j_0-1$
by Lemma \ref{nicer-lemma}.

For each $I_1$. 
the claim then follows from \eqref{suzerain}, Lemma \ref{local}, and the identity $T_{j_0-1+\m_i} = T_{j_0-1+\m_i} T_{j_0 + \m_i}$.

Now consider \eqref{d3}.  By Lemma \ref{count-again}, there are a bounded
number of functions $\phi_{I,j}^l$ with $j\le j_0$ which have support
near $3I_0$. Also, since $I_0$ is near the boundary of $\tilde{E}_{j_0}$,
we have $j\ge j_0-2$ for all these functions $\phi_{I,j}^l$.
Thus it suffices to show
$$ 
|c^l_{I,j}| 2^{\K (-j-\m_i)/2} \leq C |I_0|^{1/2} \size^*_i(T)$$
for each of the coefficients $c^l_{I,j}$ involved.

From \eqref{big-bound} and Cauchy-Schwarz we see that
$$ |c^l_{I,j}| \leq C 2^{\K(j + \m_i)/2} \| \tilde \chi_{I_0}^{10}  S_{j + \m_i} f_i \|_2 \ \ \ .$$
Since $j_0-2\le j\le j_0$ we may replace $I_0$ on the right
hand side by the $\K$-dyadic interval of length $2^{-\K j}$
which contains $I_0$.
The desired estimate follows now from Lemma \ref{get-tile-lemma} and
\eqref{unsigned-est-sup}.

The treatment of \eqref{d4} is similar to \eqref{d3}, which concludes the discussion
of the contributions \eqref{d1}, \eqref{d2}, \eqref{d3}, \eqref{d4} which come from outside 
$\tilde E_{j_0}$. 

We turn to \eqref{d5}. As with the corresponding treatment of \eqref{d5} in \eqref{mess-1}, we shall extract a gain by interacting the smoothing of $\tilde \pi_{j_0}$ with the moment condition
\eqref{mean}.

By \eqref{muj-def} and the triangle inequality\footnote{The use of the triangle inequality is a little inefficient, costing a factor of $2^{-\K(j-j_0)/2}$ or so, because it does not exploit orthogonality in physical space.  However the mean zero condition gives us a gain of $2^{-\K(j-j_0)}$, so we still end up with the improvement of $2^{-(j-j_0)/100}$ at the end.}
it suffices to show that
$$ \| \tilde \pi_{j_0} (H^l_I S_{j+\m_i} f_i - c^l_{I,j} \phi^l_{I,j}) \|_{L^2(I_0)}
\leq C |I_0|^{-1/2} \size^*_i(T) \int \tilde \chi_{I_0}^2 2^{-|j-j_0|/100} (1 + 2^{\K j} |x - x^l_I|)^{-100}$$
for all $j_0 < j $ and $I \in \Omega_j$.  Evaluating the integral, this becomes
$$ \| \tilde \pi_{j_0} (H^l_I S_{j+\m_i} f_i - c^l_{I,j} \phi^l_{I,j}) \|_{L^2(I_0)}
\leq C \size^*_i(T) 2^{\K j_0/2} 2^{-\K j} 2^{-(j-j_0)/100}
\tilde \chi_{I_0}(x^l_I)^2.$$

Fix $j$, $I$.  Observe from Fourier support considerations that
$$ \tilde \pi_{j_0} ((T_{j+\m_i-1} H^l_I) S_{j+\m_i} f_i) = 0;$$
in particular, $(T_{j+\m_i-1} H^l_I) S_{j+\m_i} f_i$ has mean zero.  It thus suffices to show that
$$ \| \tilde \pi_{j_0} F_{j,I,i} \|_{L^2(I_0)}
\leq C \size^*_i(T) 2^{\K j_0/2} 2^{-\K j} 2^{-(j-j_0)/100}
\tilde \chi_{I_0}(x^l_I)^2$$
where
$$ F_{j,I,i} := [(1 - T_{j+\m_i-1})H^l_I] S_{j+\m_i} f_i - c^l_{I,j} \phi^l_{I,j}.$$
From the construction of $c^l_{I_j}$ we see that $F_{j,I,i}$ has mean zero, and thus has a primitive $\nabla^{-1} F_{j,I,i}$ which goes to zero at $\pm \infty$.  We thus may write
$$ \| \tilde \pi_{j_0} F_{j,I,i} \|_{L^2(I_0)} = 2^{\K(j_0 + \m_i)}
\| 2^{-\K(j_0 + \m_i)} \nabla \tilde \pi_{j_0} (\nabla^{-1} F_{j,I,i}) \|_{L^2(I_0)}.$$
The multiplier $2^{-\K(j_0 + \m_i)} \nabla \tilde \pi_{j_0}$ is bounded in $L^2$ and is essentially local at scale $2^{-\K (j_0 + \m_i)}$, and hence at scale $2^{-\K j_0}$.  Thus we may estimate the previous by
\be{whoosh}
2^{\K(j_0 + \m_i)}
\| \tilde \chi_{I_0}^{100}  \nabla^{-1} F_{j,I,i} \|_2.
\end{equation}
We now claim the pointwise bound
\be{f-point}
|F_{j,I,i}(x)| \leq C \size^*_i(T) 2^{\K \m_i/2} 
(1 + 2^{\K(j+\m_i)} |x-x^l_I|)^{-50}.
\end{equation}
Assuming \eqref{f-point} for the moment, we may integrate it (using the mean zero condition to integrate from either $+\infty$ or $-\infty$, depending on which one gives the more favorable estimate) to obtain
$$
|\nabla^{-1} F_{j,I,i}(x)| \leq C \size^*_i(T) 2^{-\K(j+\m_i)} 
2^{\K \m_i/2} (1 + 2^{\K(j+\m_i)} |x-x^l_I|)^{-49}.$$
We can thus estimate \eqref{whoosh} by
$$ C 2^{\K(j_0 + \m_i)} \tilde \chi_{I_0}(x^l_I)^{100}
\size^*_i(T) 2^{-\K(j+\m_i)} 2^{\K \m_i/2} 2^{-\K(j+\m_i)/2}\ \ \ ,$$
which is acceptable.

It remains to show \eqref{f-point}.  The contribution of $c^l_{I,j} \phi^l_{I,j}$ is acceptable from \eqref{clb-2}, so it remains to consider $[(1 - T_{j+\m_i-1})H^l_I] S_{j+\m_i} f_i$.  From repeated integration by parts we have the pointwise estimate
$$ |(1 - T_{j+\m_i-1})H^l_I(x)| \leq C (1 + 2^{\K(j + \m_i)} |x-x^l_I|)^{-100}$$
so it suffices to show that
$$
|S_{j+\m_i} f_i(x)| \leq C 2^{\K \m_i/2} (1 + 2^{\K(j + \m_i)} |x-x^l_I|)^{50}.$$
Let $I'$ be a dyadic interval of length $2^{-\K j}$ which 
is adjacent to the left endpoint of $I$.
By Lemma \ref{get-tile-lemma} we can 
find a multi-tile $\pv' \in T$ 
with $I_{\pv}\subseteq 10I'$.
From \eqref{unsigned-est-sup} we thus have
$$ 
\| \tilde \chi_{I'}^{10} S_{j+\m_i} f_i \|_2 \leq C 2^{-\K j/2} \size^*_i(T).$$
From Lemma \ref{bernstein} we thus have
$$ 
\| \tilde \chi_{I'}^{10} S_{j+\m_i} f_i \|_\infty \leq C 2^{\K \m_i/2} \size^*_i(T)$$
and the desired bound follows.  This completes the treatment of \eqref{d5}.

The treatment of \eqref{d6} is similar to \eqref{d5}.
This completes the proof of \eqref{mess-2}, and therefore of \eqref{error-1}. 
The proof of Proposition \ref{proj} is thus (finally!) complete.
\endprf

\section{Deducing Theorem \ref{main} from Proposition \ref{tree-est-prop}}\label{deduce-sec}

In this section we state standard Propositions which will allow us to deduce \eqref{total-tile} and hence Theorem \ref{main} from Proposition 
\ref{tree-est-prop}.

The idea is to break the multi-tile set $\Pv_1$ into trees $T$, such that one has control on the $i$-sizes $\size^*_i(T)$ and on the total tree width $\sum_T |I_T|$.  This will be accomplished using Proposition \ref{initial} and the 
 counting function estimate of Proposition \ref{selection} on trees of a given size.

The selection of the trees is done by a greedy selection process, which will be
defined in various steps. We need the following definition:

\begin{definition}\label{convex-def}
Call a tree convex, if it is a selected tree in a greedy selection
process. In particular, convex trees satisfy Lemma \ref{simple-nesting}.
Call a subset $\Pv\subseteq \Pv_1$ convex, if it is of the  form
$\Pv_1\setminus (T_1\cup \dots\cup T_k)$ where $T_1,\dots, T_k$ are the
selected trees of a greedy selection process.
\end{definition}

\begin{proposition}\label{selection}  Let $1 \leq i \leq n$, $m \in \Z$, and suppose that $\Pv$ is a convex collection of multi-tiles such that
\be{m-bound} \size^*_i(\Pv) < \|f_i\|_2 2^{m/2}.
\end{equation}
Then there exists a collection $\T$ of distinct convex trees in $\Pv$ such that
\be{t-size}
 \sum_{T \in \T} |I_T| \leq C 2^{-m}
\end{equation}
and the remainder set $\Pv' := \Pv - \biguplus_{T \in \T} T$ is convex and satisfies
\be{remainder} 
\size^*_i(\Pv') < \|f_i\|_2 2^{(m-1)/2}.
\end{equation}
\end{proposition}

We prove this Proposition in Section \ref{select-sec}; it is the main step
in our tree selection alghorithm.  

We now aim to prove \eqref{total-tile}. By varying the $p_i$ slightly and
using Marcinkiewicz interpolation \cite{janson} it suffices to prove this estimate
under the assumption that $f_i=\chi_{E_i}$ are characteristic functions.

Starting with $m$ large and working downward, applying Proposition 
\ref{selection} for each $1 \leq i \leq n$ for each $m$, we obtain

\begin{corollary}\label{select-corollary}  For every integer $m$ there exists a collection $\T_m$ of distinct convex trees in $\Pv_1$ such that we have the size estimate
\be{size-est}
\size^*_i(\T_m) < |E_i|^{1/2} 2^{m/2}
\end{equation}
for all $1 \leq i \leq n$ and $m \in \Z$, the total tree width estimate
\be{tree-count}
\sum_{T \in \T_m} |I_T| \leq C 2^{-m}
\end{equation}
for all $m \in \Z$, and the partitioning
\be{partition}
\Pv_1 = \Pv_2\cup \biguplus_{m \in \Z} \biguplus_{T \in \T_m} T.
\end{equation}
where $\Pv_2$ is a subset of $\Pv_1$ with $\size^*_i(\Pv_2)=0$ for all $1\le i\le n$.
\end{corollary}

In the $i=n$ case we apply Proposition \ref{initial}, \eqref{size-est}, and the fact that $f_n$ is a characteristic function to obtain
\be{size-max}
\size^*_n(\T_m) \leq C \min(2^{m/2} |E_n|^{1/2}, 1)
\end{equation}
This is in fact true for all $1 \leq i \leq n$, but we shall only exploit it
for $i=n$.

Applying \eqref{partition} we may estimate the left-hand side of \eqref{total-tile} by
$$
\sum_{m \in \Z} \sum_{T \in \T_m} |\sum_{\pv \in T}
\int \chi_{I_\pv,j_\pv} \prod_{i=1}^n \pi_{\omega_{P_i}} f_i|.$$
(Observe that the set $\Pv_2$ gives no contribution, e.g. by
an appropriate application of Proposition {tree-est-prop}.)

We may apply Proposition \ref{tree-est-prop} with $\theta_i := 2/p_i$ for $1 \leq i \leq n-1$, and estimate the previous expression by
$$
C \sum_{m \in \Z} \sum_{T \in \T_m}
|I_T| (\prod_{i=1}^{n-1} \size^*_i(T)^{2/p_i}) \size^*_n(T).$$
By \eqref{size-est}, \eqref{size-max} and then \eqref{tree-count}, we may estimate this by
$$
C \sum_{m \in \Z} 2^{-m}
(\prod_{i=1}^{n-1} (2^{m/2} |E_i|^{1/2})^{2/p_i})
\min(2^{m/2} |E_n|^{1/2},1).$$
By \eqref{scaling} this simplifies to
$$
C 
(\prod_{i=1}^{n-1} |E_i|^{1/p_i})
\sum_{m \in \Z}
\min(2^{m/2 - m/p_n} |E_n|^{1/2},2^{-m/p_n}).$$
Performing the $m$ summation we obtain the desired estimate
$$
|\sum_{\pv \in \Pv_1}
\int \chi_{I_\pv,j_\pv} \prod_{i=1}^n \pi_{\omega_{P_i}} f_i| \leq C 
\prod_{i=1}^n |E_i|^{1/{p_i}}.
$$
and conclude Theorem \ref{main}.

It remains only to prove Proposition \ref{selection}.

\section{Proof of Proposition \ref{selection}}\label{select-sec}

Fix $1 \leq i \leq n$, $f_i$ and $\Pv$.  We shall need to split our notion of size $\size^*_i(\Pv)$ into upper and lower components.

For any $\xi_i$ and any sign $\pm$, let $H^\pm_{\xi_i}$ denote the Riesz projection to the half-line $\{ \xi \in \R: \pm(\xi - \xi_i) > 0\}$ in frequency space.  This Riesz projection is a linear combination of a modulated Hilbert transform and the identity.  

Observe that
$$ \|f_i \|_{P_i, \xi_i} \sim \|H^+_{\xi_i} f_i \|_{P_i, \xi_i} + \|H^-_{\xi_i} f_i \|_{P_i, \xi_i}$$
for any tile $P_i$ and any $\xi_i$.
Thus if we define
\be{size-signed-def}
\size_{i,\pm}(T) := (\frac{1}{|I_T|} \sum_{\pv \in T} \| H^\pm_{(\xi_T)_i} f_i\|_{P_i,(\xi_T)_i}^2)^{1/2}+ 
|I_T|^{-\frac 12} 
\sup_{m_{i,T}}\|\tilde{\chi}_{I_T}^{10}T_{m_{i,T}}(H^\pm_{(\xi_T)_i}f_i)\|_2,
\end{equation}
where $m_{i,T}$ is a multiplier in the range defined in Definition \ref{size-def},
and
$$
\size^*_{i,\pm}(\Pv) := \sup_{(T,\xi,I): T \subseteq \Pv} \size_{i,\pm}(T),
$$
then we have
$$ \size^*_i(\Pv) \sim \size^*_{i,+}(\Pv) + \size^*_{i,-}(\Pv).
$$
Proposition \ref{selection} will then follow from a finite number of applications of

\begin{proposition}\label{selection-alt} Let $\Pv$ be a convex collection
of multi-tiles. Let $\pm$ be a sign, and let $m \in \Z$ be such that
\be{m-bound-alt} \size^*_{i,\pm}(\Pv) < \|f_i\|_2 2^{m/2}.
\end{equation}
Then there exists a collection $\T$ of distinct convex trees in $\Pv$ such that
\be{t-size-alt}
 \sum_{T \in \T} |I_T| \leq C 2^{-m}
\end{equation}
and the remainder set $\Pv' := \Pv - \biguplus_{T \in \T} T$ is convex and satisfies
\be{remainder-alt} 
\size^*_{i,\pm}(\Pv') < \|f_i\|_2 2^{(m-1)/2}.
\end{equation}
\end{proposition}

Estimate \eqref{t-size-alt} is a variant of Bessel's inequality,
expressing that the distinct trees in this proposition correspond
to almost orthogonal components of $f_i$.

\begin{proof}
We shall prove Proposition \ref{selection-alt} for the sign +; the other sign follows by applying the frequency reflection $\xi \to -\xi$ and conjugating $f_i$.

The idea is to remove maximal 
trees in a greedy selection process 
from $\Pv$ until \eqref{remainder-alt} is obeyed for the remainder set.
This procedure shall be given by iteration.  If \eqref{remainder-alt} holds,
we terminate the iteration.  If \eqref{remainder-alt} does not hold, then there exists a tree for which
\be{big-tree}
\size_{i,+}(T) \geq \|f_i\|_2 2^{(m-1)/2}.
\end{equation}
Since the set of all possible trees $(T,\xi,I)$ obeying \eqref{big-tree} is compact, we may select $T$ so that $(\xi_T)_i$ is maximal. By retaining the top data
but adding further multi-tiles if
necessary, we may assume that $T$ is maximal in the sense of Definition
\ref{max-tree-def}.
We then add this tree $T$ to $\T$.
Then, we remove all the multi-tiles in $T$ from $\Pv$.  We then repeat this iteration until \eqref{remainder-alt} holds.

Since $\Pv$ is finite, this procedure halts in finite time and yields a collection $\T$ of mutually distinct convex trees.  Note that trees with a larger value of $\xi_T$ will be selected before trees with a smaller value of $\xi_T$.  The property \eqref{remainder-alt} holds by construction, so it only remains to show \eqref{t-size-alt}.

As usual, we shall use the $TT^*$ method to prove this orthogonality estimate.  One may think of this Lemma as a phase space version of Lemma \ref{almost}, which was set entirely in physical space, and we shall need Lemma \ref{almost} in the proof of this estimate.

Write $X := \size_{i,+}^*(\Pv)$.
Observe that
\begin{equation}\label{size-lower}
X/2\le \size_{i,+}(T) \le X
\end{equation}
for all $T\in \T$ by construction.

It suffices to prove 
\eqref{t-size-alt} separately for the set of all $T\in \T$ which satisfy
\begin{equation}\label{lacunary-size}
(\frac{1}{|I_T|} \sum_{\pv \in T} \| H^\pm_{(\xi_T)_i} f_i\|_{P_i,(\xi_T)_i}^2)^{1/2}
\ge X/4
\end{equation}
and the set of all $T\in \T$ which satisfy 
\begin{equation}\label{top-size}
|I_T|^{-\frac 12} 
\sup_{m_{i,T}}\|\tilde{\chi}_{I_T}^{10}T_{m_{i,T}}(H^\pm_{(\xi_T)_i}f_i)\|_2\ge X/4 
\end{equation}

for appropriate $m_{i,T}$.
We first consider the set of trees which satisfy \eqref{lacunary-size}. For
simplicity of notation we may assume this set is equal to $\T$.

From \eqref{size-lower}, \eqref{size-signed-def}, \eqref{fips} we may associate to each $T \in \T$ and $\pv \in T$ a multiplier $m_{P_i}$ supported on the interval
$$ \omega^+_{\pv} := \{ \xi \in 10\omega_\pv: \xi \geq (\xi_T)_i \}
$$
obeying \eqref{mpi-est} with $\xi_i = (\xi_T)_i$ such that
\be{cpv-sq}
\sum_{\pv \in T} c_{\pv}^2 \sim X^2 |I_T|,
\end{equation}
where $c_\pv$ is the non-negative quantity
$$
c_{\pv} := \| \tilde \chi_{I_\pv}^{10} T_{m_{P_i}} f_i \|_2.$$
From the signed version of \eqref{unsigned-est} we have
\be{cp-est}
c_\pv \leq C X |I_\pv|^{1/2}.
\end{equation}

Summing \eqref{cpv-sq} in $T$, we obtain
$$
\sum_{T \in \T} \sum_{\pv \in T} c_\pv c_{\pv} \sim X^2 \sum_{T\in \T}|I_T|.$$
On the other hand, from the definition of $c_{\pv}$ and duality we can find for each $\pv$ an $L^2$-normalized function $g_\pv$ such that
$$
c_\pv = \langle f_i, T_{m_{P_i}} (\tilde{\chi}_{I_\pv}^{10} g_\pv) \rangle.$$
We thus have
$$
\sum_{T \in \T} \sum_{\pv \in T} c_\pv
\langle f_i, T_{m_{P_i}} (\tilde \chi_{I_\pv}^{10} g_\pv) \rangle \sim X^2 \sum_{T \in \T} |I_T|.$$
We can write the left-hand side as
$$
\langle f_i, 
\sum_{T \in \T} \sum_{\pv \in T} c_\pv
T_{m_{P_i}} (\tilde \chi_{I_\pv}^{10} g_\pv)
\rangle.$$
By the Cauchy-Schwarz inequality we thus have
$$
X^2 \sum_{T \in \T} |I_T| \leq C \|f_i\|_2
\| \sum_{T \in \T} \sum_{\pv \in T} c_\pv
T_{m_{P_i}} (\tilde \chi_{I_\pv}^{10} g_\pv) \|_2.$$
To prove \eqref{t-size-alt} it suffices to show that
$$
\| \sum_{T \in \T} \sum_{\pv \in T} c_\pv
T_{m_{P_i}} (\tilde \chi_{I_\pv}^{10} g_\pv) \|_2 \leq C X
(\sum_{T \in \T} |I_T|)^{1/2}.
$$
It will be necessary to dyadically decompose the operator $T_{m_{P_i}}$
around the base frequency $(\xi_T)_i$.  For each $\pv \in T \in \T$, decompose
$$ T_{m_{P_i}} = \sum_{s=0}^\infty 2^{-s} T_{m_{P_i,s}}$$
where $m_{P_i,s}$ is a bump function adapted to the region
$$ \omega^+_{\pv,\xi_T,i,s} := 
[(\xi_T)_i+2^{-s}5000C_0 |\omega_{P_i}|,(\xi_T)_i+2^{2-s}5000C_0|\omega_{P_i}|]$$
and supported in $10\omega_{P_i}$.

We can then decompose 
$$ T_{m_{P_i}} (\tilde \chi_{I_\pv}^{10} g_\pv) = \sum_{s=0}^\infty 2^{-s} h_{\pv,s},$$
where 
\be{hps-def}
h_{\pv,s} := T_{m_{P_i,s}} (\tilde \chi_{I_\pv}^{10} g_\pv).
\end{equation}
Thus $h_{\pv,s}$ has Fourier support in $\omega^+_{\pv,\xi_T,i,s}$, 
is bounded in $L^2$, and is rather weakly localized in physical space near $I_\pv$.  For $s \leq M_i$ the multiplier $T_{m_{P_i,s}}$ has good spatial localization properties, but for $s > M_i$ the multiplier $T_{m_{P_i,s}}$ begins to spread $h_{\pv,s}$ along a wider interval than $I_\pv$.  However in the case $s > M_i$ we have the easily verified pointwise estimate
\be{hps-large}
|h_{\pv,s}(x)| \le C 2^{s-M_i} |I_\pv|^{-1/2} \tilde \chi_{2^{s-M_i} I_\pv}(x)^5
\end{equation}
from kernel bounds on $T_{m_{P_i,s}}$.

By the triangle inequality it thus suffices to show that
$$
\| \sum_{\pv \in T \in \T}  c_\pv
h_{\pv,s} \|_2 \leq C X 2^{s/100} (\sum_{T \in \T} |I_T|)^{1/2}
$$
for all $s \geq 0$.

Fix $s$.  We square this as
$$
\sum_{\pv \in T \in \T} \sum_{\pv' \in T' \in \T}
c_\pv c_{\pv'} \langle h_{\pv, s}, h_{\pv', s} \rangle
\leq C X^2 2^{s/50} \sum_{T \in \T} |I_T|.$$
By symmetry it suffices to show that
\be{ch-upper}
\sum_{\pv \in T \in \T} \sum_{\pv' \in T' \in \T: |\omega_\pv| \leq |\omega_{\pv'}|}
c_\pv c_{\pv'}
|\langle h_{\pv, s}, h_{\pv', s} \rangle|
\leq C X^2 2^{s/50} \sum_{T \in \T} |I_T|.
\end{equation}

We first consider the contribution when $|\omega_\pv| = |\omega_{\pv'}|$.  It suffices to show that
$$
\sum_{\pv \in T \in \T, \pv' \in T' \in \T: |\omega_{\pv}| = |\omega_{\pv'}| = 2^{\K m}}
c_\pv c_{\pv'}
|\langle h_{\pv, s}, h_{\pv', s} \rangle|
\leq C 2^{s/50} \sum_{\pv \in T \in \T: |\omega_\pv|=2^{\K m}} c_\pv^2$$
for all integers $m$, since the claim then follows by summing in $m$ and applying \eqref{cpv-sq}.

Fix $m$.  By Schur's test (i.e. estimating $c_\pv c_{\pv'} \leq \frac{1}{2}(c_\pv^2 + c_{\pv'}^2)$) and symmetry it suffices to show that
$$
\sum_{\pv' \in T' \in \T: |\omega_{\pv'}| = |\omega_\pv|}
|\langle h_{\pv, s}, h_{\pv', s} \rangle|
\leq C 2^{s/50}$$
for all $\pv \in T \in \T$.  

Fix $\pv$.  From  \eqref{sparse-single-scale}
we conclude that $\langle h_{\pv, s}, h_{\pv', s} \rangle\neq 0$ 
implies $\pv=\pv'$.
For those values one has a bound of
$$ |\langle h_{\pv, s}, h_{\pv', s} \rangle| \leq C
2^{-\K (s-M_i)_+} (1 + 2^{\K (j_\pv + (s-M_i)_+)} \dist(I_\pv, I_{\pv'}))^{-2},$$
as can be easily verified from \eqref{hps-def} when $s \leq M_i$ and \eqref{hps-large} when $s > M_i$.  The claim then follows by summing.

It remains to control the contribution when $|\omega_\pv| < |\omega_{\pv'}|$.   It suffices to show that
\be{rotd}
\sum_{T' \in \T} \sum_{\pv' \in T': |I_{\pv'}| < |I_\pv|}
|I_\pv|^{1/2} |I_{\pv'}|^{1/2} 
|\langle h_{\pv, s}, h_{\pv', s} \rangle|
\leq
C |I_\pv| (1 + \frac{\dist(2^{(s-M_i)_+} I_\pv, \R \backslash I_T)}{|I_\pv|})^{-2}
\end{equation}
for all $T \in T$ and $\pv \in T$, since we may then sum in $\pv$ to obtain
$$\sum_{\pv \in T} \sum_{T' \in \T} \sum_{\pv' \in T': |I_{\pv'}| < |I_\pv|}
|I_\pv|^{1/2} |I_{\pv'}|^{1/2} 
|\langle h_{\pv, s}, h_{\pv', s} \rangle|
\leq
C 2^{s/50} |I_T|,$$
and \eqref{ch-upper} follows by summing in $T$ and applying \eqref{cp-est}.

It remains to show \eqref{rotd}. Fix $\pv$, $T$.  Let $\Pv(\pv,s)$ denote the set of all $\pv'$ which make a non-zero contribution to \eqref{rotd}.  We now make the key observation that by Lemma \ref{plus-separation-lemma} we have
intervals $I_{\pv'}$ are disjoint from $I_T$. Namely, non-vanishing of 
$\langle h_{\pv, s}, h_{\pv', s} \rangle$ implies 
hypotheses \eqref{non-empty-10omega} and \eqref{non-empty-omega}.
Moreover, by similar reasoning, if $\pv'$ and $\pv''$ are in
$\Pv(\pv,s)$ and $\pv'\neq \pv''$, then $I_{\pv'}$ and $I_{\pv''}$
are disjoint, as one can see from \eqref{sparse-single-scale} if
$|I_{\pv}|=|I_{\pv'}|$ and from (a slight variant of) Lemma 
\ref{plus-separation-lemma} if $|I_{\pv'}|\neq |I_{\pv''}|$.
For future reference we summarize:
\begin{observation}\label{key-obs}
The intervals $I_{P'}$ are disjoint from each other and from $I_T$.
\end{observation}

We first verify \eqref{rotd}
in the case $s > M_i$.  In this case we see from \eqref{hps-large} and a calculation that
$$
|\langle h_{\pv, s}, h_{\pv', s} \rangle|
\leq C |I_{\pv'}|^{-1/2} |I_\pv|^{-1/2} 2^{-(s-M_i)} 
\int_{I_{\pv'}} \tilde \chi_{2^{s-M_i} I_\pv}^3,$$
so by Observation \ref{key-obs} the inequality \eqref{rotd}
reduces to
$$ \int_{\R \backslash I_T}
2^{-(s-M_i)} \tilde \chi_{2^{s-M_i} I_\pv}^3
\leq
C |I_\pv| (1 + \frac{\dist(2^{(s-M_i)} I_\pv, \R \backslash I_T)}{|I_\pv|})^{-2}$$
which is easily verified.

Now consider the case $s \leq M_i$.  By \eqref{hps-def} we have
$$
|\langle h_{\pv,s}, h_{\pv',s} \rangle |
\leq \langle F, \tilde \chi_{I_{\pv}}^{10} \tilde \chi_{I_{\pv'}}^{10} |g_{\pv'}| \rangle$$
where
$$ F :=  
\tilde \chi_{I_{\pv}}^{-10}
\sup_{\pv': |I_{\pv'}| \leq |I_\pv|}
|T_{m_{P'_i,s}}^* T_{m_{P_i,s}} (\tilde \chi_{I_\pv}^{10} g_\pv)|.$$
From the decay of the kernel of $T_{m_{P'_i,s}} T_{m_{P_i,s}}$ we can control $F$ pointwise by the Hardy-Littlewood maximal function:
$$ F(x) \leq C M g_\pv(x).$$

From the previous \eqref{rotd} reduces to
\be{big}
\langle M g_\pv, \tilde \chi_{I_\pv}^{10}
\sum_{\pv' \in \Pv(\pv,s) }
|I_{\pv'}|^{1/2} \tilde \chi_{I_{\pv'}}^{10} |g_{\pv'}| \rangle
\leq
C |I_\pv|^{1/2} (1 + \frac{\dist(I_\pv, \R \backslash I_T)}{|I_\pv|})^{-2},
\end{equation}
so by Cauchy-Schwarz and the Hardy-Littlewood maximal  inequality it suffices to show that
$$\| \tilde \chi_{I_\pv}^{10}
\sum_{\pv' \in \Pv(\pv,s) }
|I_{\pv'}|^{1/2} \tilde \chi_{I_{\pv'}}^{10} |g_{\pv'}| 
\|_2 
\leq
C |I_\pv|^{1/2} (1 + \frac{\dist(I_\pv, \R \backslash I_T)}{|I_\pv|})^{-2}.$$

Let us first consider the portion of the $L^2$ norm in the region
$$ \Omega := \{ x \in I_T: \dist(x, \R \backslash I_T) \geq |I_\pv| + \frac{1}{2} \dist(I_\pv, \R \backslash I_T) \}.$$
In this region we have from the $L^2$ boundedness of $g_{\pv'}$ and the decay of $\tilde \chi_{I_{\pv'}}^{10}$ that
$$ \| |I_{\pv'}|^{1/2} \tilde \chi_{I_{\pv'}}^{10} |g_{\pv'}| \|_{L^2(\Omega)}
\leq C |I_\pv|^{-1/2} \int_{I_{\pv'}} (1 + \frac{\dist(x, \Omega)}{|I_\pv|})^{-5}\ dx$$
(in fact one can get much better bounds than this, especially if $|I_{\pv'}| \ll |I_\pv|$) so by Cauchy-Schwarz and the above key observation (pairwise
disjointness of $I_{\pv'}$ and disjointness from $I_T$) again
this contribution to \eqref{big} is bounded by
$$ C |I_\pv|^{-1/2} \int_{\R \backslash I_T} (1 + \frac{\dist(x, \Omega)}{|I_\pv|})^{-5}\ dx$$
which is acceptable.

It remains to estimate the contribution outside of $\Omega$.  
Since we have
$$ \|\tilde \chi_{I_{\pv}}^{5} \chi_{\R \backslash \Omega}\|_\infty
\leq C 
(1 + \frac{\dist(I_\pv, \R \backslash I_T)}{|I_\pv|})^{-2},$$
this follows simply from Lemma \ref{almost-useful}.

This completes the proof of \eqref{t-size-alt} for the trees which satisfy
\eqref{lacunary-size}.

Now we consider the set of trees which satisfy \eqref{top-size}. For
simplicity of notation we may again assume this set is equal to $\T$.
The proof is a reprise of the previous case.

We may associate to each $T \in \T$ a multiplier $m_{i,T}$ supported on the interval
$$ \omega^+_{i,T   } := \{ \xi \in 10\omega_{i,T}: \xi \geq (\xi_T)_i \}
$$
obeying \eqref{mpi-est-tree} such that
\be{cpv-sq-tree}
 c_{T}^2 \sim X^2 |I_T|
\end{equation}
where $c_T$ is the non-negative quantity
$$
c_{T} := \| \tilde \chi_{I_T}^{10} T_{m_{i,T}} f_i \|_2.$$

By duality we can find for each $T\in \T$ an $L^2$-normalized function $g_T$ such that
$$
c_T = \langle f_i, T_{m_{i,T}} (\tilde{\chi}_{I_T}^{10} g_T) \rangle.$$
We thus have
$$
\sum_{T \in \T} c_T
\langle f_i, T_{m_{i,T}} (\tilde \chi_{I_T}^{10} g_T) \rangle \sim X^2 
\sum_{T \in \T} |I_T|.$$
By the Cauchy-Schwarz inequality we have
$$
X^2 \sum_{T \in \T} |I_T| \leq C \|f_i\|_2
\| \sum_{T \in \T} c_T
T_{m_{i,T}} (\tilde \chi_{I_T}^{10} g_T) \|_2.$$
To prove \eqref{t-size-alt} it thus suffices to show that
$$
\| \sum_{T \in \T} c_T
T_{m_{i,T}} (\tilde \chi_{I_T}^{10} g_T) \|_2 \leq C X
(\sum_{T \in \T} |I_T|)^{1/2}.
$$
We dyadically decompose the operator $T_{m_{i,T}}$
around the base frequency $(\xi_T)_i$:

$$ T_{m_{i,T}} = \sum_{s=0}^\infty 2^{-s} T_{m_{i,T,s}}$$
where $m_{i,T,s}$ is a bump function adapted to the region
$$ \omega^+_{i,T,s} := 
[(\xi_T)_i+2^{-s}10|\omega_{i,T}|,(\xi_T)_i+2^{2-s}10|\omega_{i,T}|]$$
and supported in $10\omega_{i,T}$.

Define
\be{hps-def-tree}
h_{T,s} := T_{m_{i,T,s}} (\tilde \chi_{I_T}^{10} g_T).
\end{equation}

By the triangle inequality it thus suffices to show that
$$
\| \sum_{T \in \T}  c_T
h_{T,s} \|_2 \leq C X 2^{s/2} (\sum_{T \in \T} |I_T|)^{1/2}
$$
for all $s \geq 0$. (While a better power of $s$ can be achieved, 
we shall not be ambitious here to do so.)
Fix $s$.  By squaring and using symmetry it suffices to show that
$$
\sum_{T \in \T} \sum_{T' \in \T: |\omega_{i,T}| \leq |\omega_{i,T'}|}
c_T c_{T'}
|\langle h_{T, s}, h_{T', s} \rangle|
\leq C X^2 2^{s} \sum_{T \in \T} |I_T|.
$$

Using \eqref{cpv-sq-tree} we see that it suffices to prove
for each $T\in \T$
$$
\sum_{T' \in \T: |\omega_{i,T}| \leq |\omega_{i,T'}|}
|I_{T'}|^{1/2}
|\langle h_{T, s}, h_{T', s} \rangle|
\leq C 2^{s} |I_T|^{\frac 12} .
$$

Fix $T$.
Observe
$$ \langle h_{T, s}, h_{T', s} \rangle 
=\langle T^*_{m_{i,T',s}}T_{m_{i,T,s}}\tilde{\chi}^{10}_{I_T} g_T, 
\tilde{\chi}^{10}_{I_{T'}} g_{T'} \rangle$$
and the pointwise bound
$$|T^*_{m_{i,T',s}}T_{m_{i,T,s}}\tilde{\chi}^{10}_{I_T} g_T (x)| 
\le C \tilde{\chi}^{5}_{2^{((s-M_i)_+)}I_T}(x)Mg_T(x) $$
with the Hardy Littlewood maximal function $Mg_T$. The latter
followed from
the kernel bound
$$|K(x,y)|\le C |I_T|^{-1}2^{-(s-M_i)}(1+|I_T|^{-1}2^{-(s-M_i)}|x-y|)^{-100}$$
for the kernel $K$ of the operator
$T^*_{m_{i,T',s}}T_{m_{i,T,s}}$.

Let $\T(T)$  be the set of all $T'\in T$ such that
$|\omega_{i,T}| \leq |\omega_{i,T'}|$ and
$\langle h_{T, s}, h_{T', s} \rangle\neq 0$.
By the above and the Hardy Littlewood maximal
theorem it suffices to show

$$
\|\sum_{T' \in \T(T) }
|I_{T'}|^{1/2}
\tilde{\chi}^5_{2^{((s-M_i)_+)}I_T}
\tilde{\chi}^{10}_{I_{T'}} g_{T'}\|_2
\leq C 2^{((s-M_i)_+)}|I_T|^{1/2}$$
This however follows from Lemma \ref{almost-useful} provided
we can show $\T(T)$ can be split into two subsets,
each of which has the property that the intervals $I_{T'}$
with $T'$ in the subset are pairwise disjoint. This
however follows from Lemma \ref{tree-separation-lemma}.

This completes the proof of \eqref{t-size-alt} for the trees which satisfy
\eqref{top-size}.

\end{proof}

\end{document}